%% file: HastirWinkinDochain_FunnelControl.tex
\documentclass[preprint,10pt,twocolumn,5p,authoryear]{elsarticle}

\usepackage{amsmath,amssymb,amsfonts}
\usepackage{dsfont}
\usepackage{mathrsfs}

\usepackage{xcolor}

\usepackage{bbm}

\usepackage{tikz}
\usetikzlibrary{positioning,arrows}

\tikzset{%
    block/.style={draw, fill=white, rectangle, 
            minimum height=2em, minimum width=3em},
    input/.style={inner sep=0pt},       
    output/.style={inner sep=0pt},      
    sum/.style = {draw, fill=white, circle, minimum size=2mm, node distance=1.5cm, inner sep=0pt},
    pinstyle/.style = {pin edge={to-,thin,black}}
}
 
\newtheorem{thm}{Theorem}[section]

\newtheorem{assum}{Assumption}[section]
\newtheorem{prop}{Proposition}[section]

\newtheorem{rem}{Remark}[section]

\newproof{pf}{\textbf{Proof}}
\renewcommand*{\qed}{\hfill\ensuremath{\blacksquare}}
\newcommand{\vertiii}[1]{{\vert\kern-0.25ex\vert\kern-0.25ex\vert #1 
    \vert\kern-0.25ex\vert\kern-0.25ex\vert}}

\journal{Automatica}

\begin{document}

\begin{frontmatter}

\title{Adaptive output error feedback for a class of nonlinear infinite-dimensional systems} 

\author[Namur]{Anthony Hastir}\ead{anthony.hastir@unamur.be}     
\author[Namur]{Joseph J. Winkin}\ead{joseph.winkin@unamur.be}  
\author[Louvain]{Denis Dochain}\ead{denis.dochain@uclouvain.be}

\address[Namur]{University of Namur, Department of Mathematics and Namur Institute for Complex Systems (naXys), Rue de Bruxelles 61, B-5000 Namur, Belgium}   
\address[Louvain]{Universit\'e Catholique de Louvain, Institute of Information and Communication Technologies, Electronics and Applied Mathematics (ICTEAM), Avenue Georges Lemaitre 4-6, B-1348 Louvain-La-Neuve, Belgium}


\begin{abstract}       
An adaptive funnel control method is considered for the regulation of the output for a class of nonlinear infinite-dimensional systems on real Hilbert spaces. After a decomposition of the state space and some change of variables related to the Byrnes-Isidori form, it is shown that the funnel controller presented in \citep{SchwenningerFunnel} achieves the control objective under some assumptions on the nonlinear system dynamics, like well-posedness and Bounded-Input-State Bounded-Output (BISBO) stability. The theory is applied to the regulation of the temperature in a chemical plug-flow tubular reactor whose reaction kinetics are modeled by the Arrhenius nonlinearity. Furthermore a damped sine-Gordon model is shown to fit the required assumptions as well. The theoretical results are illustrated by means of numerical simulations. 
\end{abstract}

\end{frontmatter}

\section{Introduction}

An adaptive funnel controller is a model-free output error feedback controller with a time-varying gain that adapts according to the output error in order to ensure a predetermined transient behavior for the tracking error with no necessary asymptotic tracking performance. This field of control has attracted a lot of attention in the last few years. It has been thoroughly developed in \cite{Ilchmann2002} for systems with relative degree one. In particular, they have shown that the funnel control approach is feasible for linear finite-dimensional systems, infinite-dimensional linear regular systems, nonlinear finite-dimensional systems, nonlinear delay systems and systems with hysteresis. Details for nonlinear systems with relative degree one are also available in \cite{Ilchmann2005}. Few years later, funnel control for MIMO nonlinear systems with known strict relative degree has been widely developed in \cite{BergerFunnelAutomatica}. The new funnel controller they introduced involves the $r-1$ derivatives of the tracking error, where $r$ stands for the relative degree of the system.  As an illustration of the method, they consider a mass-spring system mounted on a car. Note that a new survey for funnel control for different types of systems can be found in \cite{Berger2021Bis}. Funnel control has also been considered in several aplications. For instance \cite{Ilchmann2004} developed a funnel controller to regulate the reaction temperature in chemical reactors. Later, the position of a moving water tank has been controlled via a funnel controller in \cite{Berger2020funnel}. A linearized version of the Saint-Venant Exner infinite-dimensional dynamics has been used as internal dynamics. This shows that funnel control becomes more and more attractive for systems driven by infinite-dimensional internal dynamics. In that way, this topic has been considered in \cite{SchwenningerFunnel}. The authors proved that some class of infinite-dimensional linear systems fits the required assumptions for funnel control to be feasible. Moreover, they proved that linear-infinite dimensional systems that can be written in Byrnes-Isidori form, see \cite{IlchmannByrnesIsidori}, are encompassed in that new class. More recently, funnel control has also been applied to a nonlinear infinite-dimensional reaction-diffusion equation coupled with the nonlinear Fitzhugh-Nagumo model, which represent together defibrillation processes of the human heart, see \cite{Berger2021}. Note also that funnel control has been lately coupled to model-predictive-control (Funnel MPC) for nonlinear systems with relative degree one, see \cite{Berger2021funnel}. 

The novelty in this paper consists in considering a class of nonlinear infinite-dimensional systems to which funnel control can be applied\footnote{Note that this class is restricted to systems with relative degree one.}. Based on the Byrnes-Isidori form for linear systems, we introduce a change of variables that aims at extracting the output dynamics of the system, which is assumed to be finite-dimensional. Based on this transformation, funnel control is shown to be feasible provided that the remaining part of the dynamics satisfies some BISBO stability assumption. This constitutes the main result of the paper. Moreover, a way of getting this BISBO stability condition is presented for bounded nonlinearities. We illustrate our results on two nonlinear infinite-dimensional partial differential equations (PDEs) modeling the dynamics of a plug-flow tubular reactor and a damped sine-Gordon nonlinear equation, respectively.

The paper is organized as follows: the system class for which funnel control is known to be feasible together with the control objective are presented in Section \ref{SystemClass}. In Section \ref{NonlinearDPS} a class of nonlinear controlled and observed infinite-dimensional systems with some appropriate assumptions is considered and is shown to fulfill the assumptions needed to apply funnel control. A way to get BISBO stability of the internal dynamics is presented. The application of the main results with some numerical simulations are given in Section \ref{Simus} on a chemical plug-flow tubular reactor and a damped sine-Gordon equation. Section \ref{Conclusions} is dedicated to conclusions and perspectives.

\section{System class and control objective}\label{SystemClass}
This section is dedicated to the introduction of the considered system class together with the control objective that is pursued. 

\subsection{General framework}
We shall consider a differential equation that makes the connection between the input and the output of a dynamical system whose internal dynamics are not necessarily known (model-free). This equation reads as
\begin{equation}
\left\{
\begin{array}{l}
\dot{y}(t) = N(d(t), T(y)(t)) + \Gamma(d(t),T(y)(t))u(t),\\
y(0) = y_0,
\end{array}\right.
\label{OutputEquation}
\end{equation}
where the following conditions are assumed to hold.
\begin{assum}\label{Disturb}
The disturbance $d\in L^\infty(\mathbb{R}^+,\mathbb{R})$, the nonlinear function $N$ is in $\mathcal{C}(\mathbb{R}^2,\mathbb{R})$ and the gain function $\Gamma\in\mathcal{C}(\mathbb{R}^2,\mathbb{R})$ is positive in the sense that $\Gamma(d,\varrho)>0$ for all $(d,\varrho)\in\mathbb{R}^2$.
\end{assum}

\begin{assum}\label{NonlinearMapT}
The map $T:\mathcal{C}(\mathbb{R}^+,\mathbb{R})\to L^\infty_{\text{loc}}(\mathbb{R}^+,\mathbb{R})$ is a (possibly nonlinear) operator which satisfies the following conditions:
\begin{enumerate}
\item Bounded trajectories are mapped into bounded trajectories (BIBO property), i.e. for all $k_1>0$, there exists $k_2>0$ such that for all $y\in\mathcal{C}(\mathbb{R}^+,\mathbb{R})$,
\begin{equation}
\sup_{t\in\mathbb{R}^+} \vert y(t)\vert\leq k_1 \Rightarrow \sup_{t\in\mathbb{R}^+} \vert T(y)(t)\vert\leq k_2
\label{BIBOAssum}
\end{equation}
\item The operator $T$ is causal, i.e. for any $t\in\mathbb{R}^+$ and any $y,\hat{y}\in\mathcal{C}(\mathbb{R}^+,\mathbb{R})$
\begin{equation}
y_{\vert [0,t)} = \hat{y}_{\vert [0,t)} \Rightarrow T(y)_{\vert [0,t)} = T(\hat{y})_{\vert [0,t)},
\label{Causality}
\end{equation}
where $f_{\vert I}$ denotes the restriction of a function $f$ to the interval $I$.
\item $T$ is locally Lipschitz in the sense that for all $t\in\mathbb{R}^+$ and all $y\in\mathcal{C}([0,t],\mathbb{R})$ there exist positive constants $\tau,\delta$ and $\rho$ such that for any $y_1,y_2\in\mathcal{C}(\mathbb{R}^+,\mathbb{R})$ with $y_{i_{\vert [0,t]}} = y, i=1,2$ and $\vert y_i(s)-y(t)\vert < \delta$ for all $s\in[t,t+\tau]$ and $i=1,2$ it holds that
\begin{equation}
\Vert (T(y_1) - T(y_2))_{\vert[t,t+\tau]}\Vert_\infty \leq \rho\Vert (y_1-y_2)_{\vert[t,t+\tau]}\Vert_\infty,
\label{Lipschitz}
\end{equation}
where $\Vert f_{\vert[t,t+\tau]}\Vert_\infty := \sup_{s\in[t,t+\tau]}\vert f(s)\vert$.
\end{enumerate}
\end{assum}
The class of systems governed by (\ref{OutputEquation}) with Assumptions \ref{Disturb}--\ref{NonlinearMapT} is presented in \cite[Section 1]{SchwenningerFunnel} for systems with (possible) memory and relative degree $r\in\mathbb{N}$. Here we consider systems with no memory and relative degree one, see (\ref{OutputEquation}). This class is quite general and encompasses systems with infinite-dimensional internal dynamics as shown in \cite{SchwenningerFunnel} and \cite{IlchmannByrnesIsidori} for instance. However it is still not clear which classes of distributed-parameter systems (DPS) may be written as the input-output equation (\ref{OutputEquation}). Such a compatible class is presented in Section \ref{NonlinearDPS}.

\subsection{Control objective}
The objective here consists in developing an output error feedback $u(t) = \mathcal{G}(t,e(t))$ with $e(t) = y(t) - y_{\text{ref}}(t)$ for a reference signal $y_{\text{ref}}\in W^{1,\infty}(\mathbb{R}^+,\mathbb{R})$, such that, when connected to (\ref{OutputEquation}), it results in a closed-loop system for which the error $e(t)$ evolves in a prescribed performance funnel
\begin{equation}
\mathcal{F}_\phi := \left\{(t,e)\in\mathbb{R}^+\times\mathbb{R},\, \phi(t)\vert e(t)\vert < 1\right\},
\label{FunnelSet}
\end{equation}
where the function $\phi$ is assumed to belong to
\begin{align}
\Phi &:= \left\{\phi\in\mathcal{C}(\mathbb{R}^+,\mathbb{R}),\, \phi, \dot{\phi}\in L^\infty(\mathbb{R}^+,\mathbb{R}),\right.\nonumber\\
&\hspace{1cm}\left.\phi(t)> 0, \forall
t\in\mathbb{R}^+ \text{ and } \liminf_{t\to\infty}\phi(t) > 0\right\}.
\label{BoundaryFunnel}
\end{align}
This control objective is also considered in \cite[Section 1]{SchwenningerFunnel}. As described in \cite{SchwenningerFunnel,IlchmannByrnesIsidori} or \cite{BergerFunnelAutomatica}, a controller that achieves the output tracking performance described above is expressed as
\begin{equation}
u(t) = \frac{-e(t)}{1-\phi^2(t)e^2(t)},
\label{FunnelController}
\end{equation}
with $\phi\in\Phi$ and $\phi(0)\vert e(0)\vert < 1$. The controller (\ref{FunnelController}) is called a funnel controller and can be viewed as the output error feedback $u(t) = -k(t)e(t)$ with a time-varying (adaptive) gain $k(t) = \frac{1}{1-\phi^2(t)e^2(t)}$. Let us consider the following theorem, coming from \cite{BergerFunnelAutomatica} with $r=1$, which characterizes the effectiveness of the controller (\ref{FunnelController}) in terms of existence and uniqueness of solutions of the closed-loop systems and also in terms of output tracking performance.
\begin{thm}\label{ThmFunnel}
Consider a system (\ref{OutputEquation}) with Assumptions \ref{Disturb}--\ref{NonlinearMapT}. Let $y_{\text{ref}}\in W^{1,\infty}(\mathbb{R}^+,\mathbb{R}), \phi\in\Phi$ and $y_0\in\mathbb{R}$ such that the condition $\phi(0)\vert e(0)\vert < 1$ holds true. Then the funnel controller (\ref{FunnelController}) applied to (\ref{OutputEquation}) results in a closed-loop system whose solution $y:[0,\omega)\to\mathbb{R}, \omega\in (0,\infty]$, has the following properties:
\begin{enumerate}
\item The solution is global, i.e. $\omega = \infty$;
\item The input $u:\mathbb{R}^+\to\mathbb{R}$, the gain function $k:\mathbb{R}^+\to\mathbb{R}$ and the output $y:\mathbb{R}^+\to\mathbb{R}$ are bounded;
\item The tracking error $e:\mathbb{R}^+\to\mathbb{R}$ evolves in the funnel $\mathcal{F}_\phi$ and is bounded away from the funnel boundaries in the sense that there exists $\epsilon> 0$ such that, for all $t\geq 0, \vert e(t)\vert \leq \frac{1}{\phi(t)} - \epsilon$.
\end{enumerate}
\end{thm}
This theorem will be a paramount tool in the next section, in order to prove output tracking control of the scalar output of some class of nonlinear DPS.

The main contribution here relies on the fact that we consider a class of nonlinear infinite-dimensional systems that admits an input-output differential description of the form (\ref{OutputEquation}). This consitutes also the difference with respect to e.g. \cite{SchwenningerFunnel} wherein a class of operators $T$ is introduced, which comes from systems modeled by linear infinite-dimensional internal dynamics. Our contribution enlarges the class of systems for which funnel control is feasible since, to the best of our knowledge, our class of nonlinear infinite-dimensional systems is shown to be appropriate for funnel control for the first time here. However, examples in which nonlinear infinite-dimensional systems are considered have already been studied for funnel control in \cite{Berger2021}. Our work also extends the Byrnes Isidori forms studied in \cite{IlchmannByrnesIsidori} to nonlinear infinite-dimensional systems that satsisfy some relatively standard assumptions.

\section{Identification of a class of nonlinear distributed parameter systems}\label{NonlinearDPS}

In this section a quite general class of nonlinear infinite-dimensional systems is introduced for which output error control is considered. It is shown how the abstract differential equation governing the dynamics can be transformed into (\ref{OutputEquation}) by using a change of variables as it is made for linear systems in \cite{IlchmannByrnesIsidori}. Moreover, by stating some quite easily checkable assumptions, it is proved that the transformed equation satisfies Assumptions \ref{Disturb}--\ref{NonlinearMapT} introduced in Section \ref{SystemClass}, which implies that funnel control is feasible for such a class of systems. This consitutes the main contributions of the paper.

Let $H$ be a real (separable) Hilbert space equipped with the inner product $\langle\cdot,\cdot\rangle_H$. The nonlinear systems that we consider are governed by the following controlled and perturbed abstract ordinary differential equation
\begin{equation}
\Sigma:\left\{
\begin{array}{l}
\dot{x}(t) = Ax(t) + f(x(t)) + b(u(t) + d(t)),\\
y(t) = \langle x(t),c\rangle_H,\\
x(0) = x_0\in H,
\end{array}\right.
\label{AbstractODE}
\end{equation}
for which we make the following three assumptions.
\begin{assum}\label{OpA}
The (unbounded) linear operator $A:D(A)\subset H\to H$ is the infinitesimal generator of a strongly continuous ($C_0$) semigroup of bounded linear operators on $H$.
\end{assum}

\begin{assum}\label{OpF}
The nonlinear operator $f:H\to H$ is uniformly Lipschitz continuous on $H$.% and is supposed to satisfy $\Vert f(0)\Vert \leq \sigma, \sigma > 0$.
%and satisfies $\Vert f(x)\Vert_H\leq \sigma$ for all $x\in H$ and $\sigma > 0$, independent of $x$.
\end{assum}

\begin{assum}\label{ShapeFct}
The vectors $b$ and $c$ are in $D(A)$ and $D(A^*)$, respectively, where $D(A^*)$ is the domain of the adjoint of the operator $A$. Moreover, it is assumed that $\langle b,c\rangle_H > 0$. The disturbance\footnote{By entering in the system (\ref{AbstractODE}) in this way, the disturbance $d$ may be interpreted as modeling uncertainties on the input, that could be due to unmodelled actuator dynamics, for instance.} $d\in L^\infty(\mathbb{R}^+,\mathbb{R})$.
\end{assum}

Note that the scalar functions $u$ and $y$ stand for the input and the output, respectively. According to \cite[Theorem 11.1.5]{CurtainZwartNew}, Assumptions \ref{OpA} and \ref{OpF} implies that the homogenenous\footnote{By "homogenenous" we mean uncontrolled and unperturbed.} part of (\ref{AbstractODE}) has a unique mild solution on $[0,\infty)$ which can even be a classical one provided that $x_0\in D(A)$. 
\begin{rem}\label{RemarkShapeFct}
We emphasize that, despite the fact that Assumption \ref{ShapeFct} may be seen as restrictive, if $b$ and $c$ are shape functions that do not lie in $D(A)$ and $D(A^*)$, it is always possible to approximate them, as accurately as desired, by functions in $D(A)$ and $D(A^*)$, since they are both dense subspaces of $H$.
\end{rem}
The consequences of Remark \ref{RemarkShapeFct} will be discussed later in this section.

According to \cite{IlchmannByrnesIsidori} and \cite{Byrnes1998}, Assumption \ref{ShapeFct} allows us to decompose the state space $H$ as
\begin{equation}
H := \text{span}\{c\} \oplus \{b\}^\perp =: \mathcal{O} \oplus \mathcal{I}.
\label{DecompositionH}
\end{equation}
To do so, let us consider the operator $P : H\to\mathbb{R}$ defined as
\begin{equation*}
Px = \frac{\langle x,b\rangle}{\langle c,b\rangle}.
\end{equation*}
This allows to write the linear subspaces $\mathcal{O}$ and $\mathcal{I}$ as $cPH$ and $(I - cP)H$, respectively. In what follows we shall denote the operator $I-cP$ by $P_\mathcal{I}$. Consequently, according to the decomposition (\ref{DecompositionH}), any element $x\in H$ may be written as 
\begin{equation*}
x = cPx + P_\mathcal{I}x.
\end{equation*} 
Moreover, observe that for any function $x\in H$ there holds $\langle P_\mathcal{I}x,b\rangle_H = 0$. The operator $P_\mathcal{I}$ is not an orthogonal projector on $\mathcal{I}$. Such a projector is denoted by $P^\perp:H\to\mathcal{I}$ and is defined as $P^\perp x = x - b\frac{\langle x,b\rangle_H}{\langle b,b\rangle_H}$. Now we shall introduce some operators that will be of importance in order to transform system (\ref{AbstractODE}). Let us consider the operator $U:H\to \mathbb{R}\times \mathcal{I}$ defined as
\begin{equation}
Ux = \left(\begin{matrix}
Px\\
P_\mathcal{I}x
\end{matrix}\right).
\label{U}
\end{equation}
This operator is boundedly invertible, see \cite{IlchmannByrnesIsidori}, with inverse $U^{-1}:\mathbb{R}\times \mathcal{I}\to H$ expressed as
\begin{equation}
U^{-1}\left(\begin{matrix}
\alpha\\
\eta
\end{matrix}\right) = \alpha c + \eta.
\label{Uinv}
\end{equation}
The adjoint of $U^{-1}$, denoted by $U^{-*}:H\to \mathbb{R}\times \mathcal{I}$, reads as follows
\begin{equation}
U^{-*}x = \left(\begin{matrix}
\langle x,c\rangle_H\\
P^\perp x
\end{matrix}\right).
\label{UinvAdj}
\end{equation}
Consequently, the operator $U^*:\mathbb{R}\times\mathcal{I}\to H$ is given by
\begin{equation}
U^*\left(
\begin{matrix}
\alpha\\
\eta
\end{matrix}
\right) = \alpha \frac{b}{\langle c,b\rangle_H} + \eta - \langle c,\eta\rangle_H\frac{b}{\langle c,b\rangle_H}.
\label{Uadj}
\end{equation} 
Like the operators $U$ and $U^{-1}$, the operators $U^*$ and $U^{-*}$ are linear and bounded. This allows us to perform the change of variables for system (\ref{AbstractODE}) defined by $x = U^*\xi, \xi\in\mathbb{R}\times\mathcal{I}$. By the invertibility of $U^*$, one gets
\begin{equation}
\xi = U^{-*}x = \left(\begin{matrix} \langle x,c\rangle_H\\ P^\perp x\end{matrix}\right).
\label{ChangeVar}
\end{equation}
By setting $\eta := P^\perp x$, one has that $\xi = (\begin{smallmatrix} y\\ \eta\end{smallmatrix})$. Now we take a look at the dynamics of the variable $\xi$. In the new variables (\ref{ChangeVar}), the system (\ref{AbstractODE}) reads as
\begin{equation}
\tilde{\Sigma}: \left\{
\begin{array}{l}
\dot{\xi}(t) = U^{-*}AU^*\xi(t) + U^{-*}f(U^*\xi(t)) + U^{-*}bu(t)\\
\hspace{3.5cm} + U^{-*}bd(t),\\
y(t) = (\begin{matrix} 1 & 0\end{matrix})\xi(t),\\
\xi(0) = U^{-*}x_0 =: \xi_0.
\end{array}\right.
\label{NonlinearDPSTransf}
\end{equation}
Note that the systems $\Sigma$ and $\tilde{\Sigma}$ are equivalent according to the state transformation induced by $U^{-*}$. By using the definitions of $U^{-*}$ and $U^{*}$, see (\ref{UinvAdj}) and (\ref{Uadj}), the linear operator $U^{-*}AU^*$ can be rewritten as
\begin{align*}
U^{-*}AU^*\xi = \left(\begin{smallmatrix}\langle AU^*\xi,c\rangle_H\\
P^\perp AU^*\xi\end{smallmatrix}\right) &= \left(\begin{smallmatrix}\langle \xi,UA^*c\rangle_{\mathbb{R}\times\mathcal{I}}\\
P^\perp AU^*\xi\end{smallmatrix}\right)\\
&=\left(\begin{smallmatrix}y\frac{\langle A^*c,b\rangle_H}{\langle c,b\rangle_H} + \langle\eta,P_\mathcal{I}A^*c\rangle_H\\
y\frac{P^\perp Ab}{\langle c,b\rangle_H} + P^\perp A\eta - \frac{\langle c,\eta\rangle_H}{\langle c,b\rangle_H}P^\perp Ab\end{smallmatrix}\right)\\
&=: \left(\begin{smallmatrix}P_0 & S\\
R & Q\end{smallmatrix}\right)\xi,
\end{align*} 
where the operators $P_0, S, R$ and $Q$ are defined as follows
\begin{align}
&P_0: \mathbb{R}\to\mathbb{R}, P_0y = y\frac{\langle A^*c,b\rangle_H}{\langle c,b\rangle_H},\nonumber\\
&S: \mathcal{I}\to\mathbb{R}, S\eta = \langle\eta,P_\mathcal{I}A^*c\rangle_H,\nonumber\\
&R: \mathbb{R}\to\mathcal{I}, Ry = y\frac{P^\perp Ab}{\langle c,b\rangle_H},\nonumber\\
&Q: D(Q)\subset\mathcal{I}\to\mathcal{I}, Q\eta = P^\perp A\eta - \frac{\langle c,\eta\rangle_H}{\langle c,b\rangle_H}P^\perp Ab,\label{OpQ}
\end{align}
where $D(Q) = D(A)\cap\mathcal{I}$. According to \cite{IlchmannByrnesIsidori}, the operator $Q$ is the infinitesimal generator of a $C_0-$semigroup\footnote{Without loss of generality there exist constants $M\geq 1$ and $\omega\in\mathbb{R}$ such that $\Vert T_Q(t)\Vert \leq Me^{\omega t}$.} $(T_Q(t))_{t\geq 0}$ on $\mathcal{I}$. Moreover, since the operators $P_0, S$ and $R$ are bounded, the operator $\left(\begin{smallmatrix}P_0 & S\\
R & Q\end{smallmatrix}\right)$ is still the generator of a $C_0-$semigroup, see e.g. \cite{CurtainZwartNew, Pazy, Engel}. 
%We shall now make the following assumption which will be needed in the main theorem of the paper.
%\begin{assum}
%The linear $C_0-$semigroup $(T_Q(t))_{t\geq 0}$ is assumed to be exponentially stable on $\mathcal{I}$, i.e. there exist constants $M \geq 1$ and $\omega > 0$ such that 
%\begin{equation}
%\Vert T_Q(t)\Vert\leq M e^{-\omega t}, t\geq 0.
%\label{ExpStab}
%\end{equation} 
%\end{assum}
By using (\ref{UinvAdj}) and (\ref{Uadj}), we can write the nonlinear part of (\ref{NonlinearDPSTransf}) as
\begin{align*}
U^{-*}f(U^*\xi(t)) &= \left(\begin{smallmatrix}\langle f(U^*\xi(t)),c\rangle_H\\
f(U^*\xi(t)) - \frac{\langle f(U^*\xi(t)),b\rangle_H}{\langle b,b\rangle_H}b
\end{smallmatrix}\right)\\
&=: \left(\begin{smallmatrix}\langle \tilde{f}(y(t),\eta(t)),c\rangle_H\\
\tilde{f}(y(t),\eta(t)) - \frac{\langle \tilde{f}(y(t),\eta(t)),b\rangle_H}{\langle b,b\rangle_H}b
\end{smallmatrix}\right),
\end{align*} 
where the nonlinear operator $\tilde{f}:\mathbb{R}\times\mathcal{I}\to H$ is defined as $\tilde{f}(y,\eta) = f(y\frac{b}{\langle c,b\rangle_H} + \eta - \langle c,\eta\rangle_H\frac{b}{\langle c,b\rangle_H})$. By using Assumption \ref{OpF} and the fact that $U^{-*}$ and $U^*$ are linear bounded operators, the nonlinear operator $U^{-*}f(U^*\cdot)$ is still uniformly Lipschitz continuous from $\mathbb{R}\times\mathcal{I}$ into $\mathbb{R}\times\mathcal{I}$. This entails that the homogeneous part of (\ref{NonlinearDPSTransf}) possesses a unique mild solution on $[0,\infty)$. Taking any initial condition in $\mathbb{R}\times D(Q)$ implies that this solution is classical. Now observe that the term $U^{-*}b$ is expressed as $\left(\begin{smallmatrix}\langle b,c\rangle_H & 0\end{smallmatrix}\right)^T$. From these observations it follows that the dynamics of $y$ and $\eta$ may be written as
\small{
\begin{align}
\dot{y}(t) = P_0y(t) + S\eta(t) + \langle \tilde{f}(y(t),\eta(t)),c\rangle_H + \gamma u(t)+ \gamma d(t)\label{NewDynamicsY}
\end{align}
}\normalsize
and
\small{
\begin{align}
\dot{\eta}(t) = Ry(t) + Q\eta(t) + \tilde{f}(y(t),\eta(t)) - \frac{\langle \tilde{f}(y(t),\eta(t)),b\rangle_H}{\langle b,b\rangle_H}b,\label{NewDynamicsEta}
\end{align}
}\normalsize
with initial condition $y(0) = y_0$ and $\eta(0) = \eta_0$, respectively, where $\gamma := \langle b,c\rangle_H$. This entails that (\ref{NewDynamicsY}) admits the representation (\ref{OutputEquation}), where
\begin{enumerate}
\item The gain function $\Gamma:\mathbb{R}^2\to\mathbb{R}$ is defined as $\Gamma(d,\varrho)=\gamma > 0$;
\item The well-defined nonlinear operator $T:\mathcal{C}(\mathbb{R}^+,\mathbb{R})\to L^\infty_{\text{loc}}(\mathbb{R}^+,\mathbb{R})$ has the form 
\begin{equation}
T(y)(t) = P_0y(t) + S\eta(t) + \langle \tilde{f}(y(t),\eta(t)),c\rangle_H,
\label{EquationT}
\end{equation}
where $\eta(t)$ is the mild solution of (\ref{NewDynamicsEta});
\item The function $N:\mathbb{R}^2\to \mathbb{R}$ reads as $N(d,\varrho) = \gamma d + \varrho$.
\end{enumerate}
This shows that Assumption \ref{Disturb} on system (\ref{OutputEquation}) is satisfied. It remains to show that the nonlinear operator $T$ given by (\ref{EquationT}) possesses the three properties of Assumption \ref{NonlinearMapT}. This constitutes the main result of this section. Before going into the details of this result, we shall denote by $\Sigma_{y\to\eta}$ the system which can be viewed as a system with input $y$ and output $\eta$ and whose dynamics are described by (\ref{NewDynamicsEta}). We make the following assumption on that system.

\begin{assum}\label{BIBO_Assum}
The system $\Sigma_{y\to\eta}$ whose dynamics are governed by (\ref{NewDynamicsEta}) is BISBO stable in the following sense: for all $k>0$ and all $\hat{k}>0$, there exists $\tilde{k}>0$ such that for all $y\in\mathcal{C}(\mathbb{R}^+,\mathbb{R})$ and all $\eta_0\in\mathcal{I}$,
\begin{equation}
\sup_{t\in\mathbb{R}^+}\vert y(t)\vert\leq k \text{ and } \Vert\eta_0\Vert\leq\hat{k}\Rightarrow \sup_{t\in\mathbb{R}^+} \Vert \eta(t)\Vert \leq \tilde{k}.
\label{BIBO_Y_ETA}
\end{equation}
\end{assum}

\begin{thm}
The operator $T$ defined in (\ref{EquationT}), which arises from the nonlinear system (\ref{AbstractODE}) via the change of variables (\ref{ChangeVar}), satisfies Assumption \ref{NonlinearMapT}.  
\end{thm}

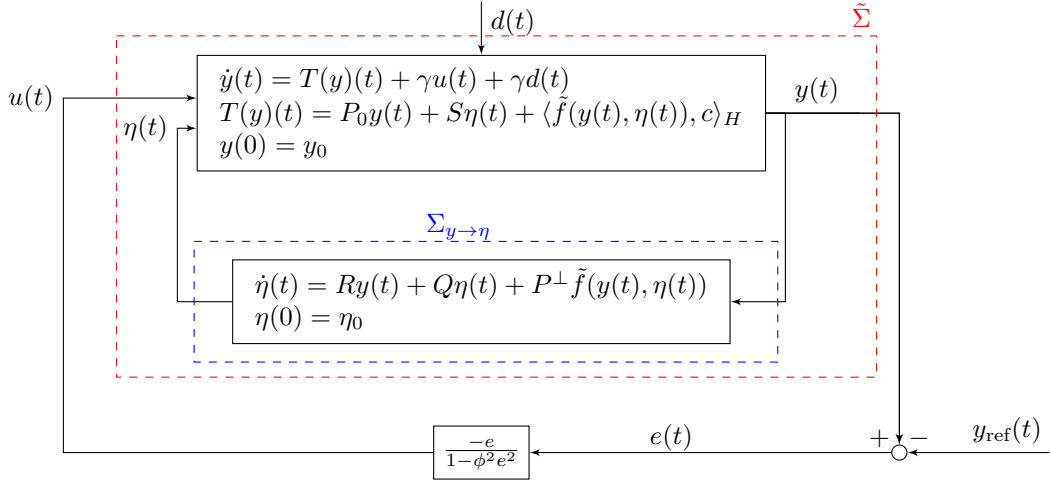
\begin{figure*}
%\begin{center}
%\includegraphics[scale=1.5,trim=0cm 1.5cm 0cm 0cm]{SchemasTikz/FunnelController.pdf}
%\caption{Interconnection of $\tilde{\Sigma}$ and the funnel controller (\ref{FunnelController}).\label{fig:closed-loop}}
%\end{center}
\begin{center}
\input{FunnelControllerCopy.tex}
\end{center}
\caption{Interconnection of $\tilde{\Sigma}$ and the funnel controller (\ref{FunnelController}).\label{fig:closed-loop}}
\end{figure*}

\begin{pf}
The proof is divided into three steps, according to the three items of Assumption \ref{NonlinearMapT}.\\
\textit{Step 1: } In order to show that $T$ maps bounded trajectories into bounded ones, let us fix $k_1 > 0$, $k_* > 0$, $y\in\mathcal{C}(\mathbb{R}^+,\mathbb{R})$ and $\eta_0\in\mathcal{I}$ such that $\sup_{t\in\mathbb{R}^+}\vert y(t)\vert\leq k_1$ and $\Vert\eta_0\Vert\leq k_*$. There exists a positive $\tilde{k}$ such that, for this $y$ and this $\eta_0$, the mild solution of (\ref{NewDynamicsEta}) with initial condition $\eta_0$ satisfies $\sup_{t\in\mathbb{R}^+} \Vert\eta(t)\Vert\leq \tilde{k}$, according to Assumption \ref{BIBO_Assum}.
Thanks to the expression (\ref{EquationT}) of $T$, the boundedness of the operator $S$ and the Cauchy-Schwarz inequality, one may write that 
\begin{align*}
\vert T(y)(t)\vert\leq \vert P_0\vert &\vert y(t)\vert + \Vert S\Vert_{\mathcal{L}(\mathcal{I},\mathbb{R})}\Vert\eta(t)\Vert\\
&\hspace{2cm} + \Vert\tilde{f}(y(t),\eta(t))\Vert_H\Vert c\Vert_H.
\end{align*}
Assumption \ref{OpF} allows us to write
\begin{align}
\vert T(y)&(t)\vert\leq \vert P_0\vert \vert y(t)\vert + \Vert S\Vert_{\mathcal{L}(\mathcal{I},\mathbb{R})}\Vert\eta(t)\Vert\nonumber\\
&\hspace{0cm} + \Vert \tilde{f}(y(t),\eta(t)) - \tilde{f}(0,0)\Vert_H\Vert c\Vert_H + \Vert\tilde{f}(0,0)\Vert_H\Vert c\Vert_H\nonumber\\
&\leq \vert P_0\vert \vert y(t)\vert + \Vert S\Vert_{\mathcal{L}(\mathcal{I},\mathbb{R})}\Vert\eta(t)\Vert\nonumber\\
&\hspace{0cm} + (l_1\vert y(t)\vert + l_2\Vert\eta(t)\Vert)\Vert c\Vert_H + \sigma\Vert c\Vert_H,\label{LipschitzConstants}
\end{align}
where $l_1>0$ and $l_2>0$ denote the Lipschitz constants of the operator $\tilde{f}$ associated with $y$ and $\eta$, respectively, and where the positive constant $\sigma$ is such that\footnote{This is valid since the nonlinear operator $f$ maps the whole space $H$ into itself, meaning that any point in $H$ has a well-defined image by $f$ in the $H-$norm.} $\Vert f(0)\Vert_H\leq \sigma$. Consequently,
\begin{align*}
\sup_{t\in\mathbb{R}^+}\vert T(y)(t)\vert \leq \vert P_0\vert k_1 + \Vert S\Vert_{\mathcal{L}(\mathcal{I},\mathbb{R})}\tilde{k} + \tilde{\sigma}\Vert c\Vert_H =: k_2,
\end{align*}
where $\tilde{\sigma} = \sigma + l_1k_1 + l_2\tilde{k}$, which proves that $T$ satisfies the BIBO condition required in Assumption \ref{NonlinearMapT}.\\
\textit{Step 2: } The causality can be easily established by noting that, for a fixed $y\in\mathcal{C}(\mathbb{R}^+,\mathbb{R})$ the corresponding mild solution of (\ref{NewDynamicsEta}) is unique. This entails that for $y,\hat{y}\in\mathcal{C}(\mathbb{R}^+,\mathbb{R})$ such that $y_{\vert [0,t)} = \hat{y}_{\vert [0,t)}$, the corresponding mild solutions of (\ref{NewDynamicsEta}), denoted by $\eta$ and $\hat{\eta}$, respectively, satisfy $\eta_{\vert [0,t)} = \hat{\eta}_{\vert [0,t)}$. In view of the expression (\ref{EquationT}) of $T$, it follows that (\ref{Causality}) holds.\\
\textit{Step 3: } For the local Lipschitz continuity, let us consider $t\geq 0$ and $y\in\mathcal{C}([0,t],\mathbb{R})$. Now let us take $y_1,y_2\in\mathcal{C}(\mathbb{R}^+,\mathbb{R})$ such that $y_i$ coincides with $y$ up to time $t$ for $i=1,2$. The mild solutions of (\ref{NewDynamicsEta}) with input $y_i$ and starting at time $t$ are given by 
\begin{align*}
\eta_i(\tilde{t}) &= T_Q(\tilde{t}-t)\eta_{i,t} + \int_t^{\tilde{t}}T_Q(\tilde{t}-s)Ry_i(s)ds\\
&\hspace{2cm}+\int_t^{\tilde{t}}T_Q(\tilde{t}-s)P^\perp\tilde{f}(y_i(s),\eta_i(s))ds
\end{align*}
for any $\tilde{t}\in [t,t+\tau]$ with $\tau$ being an arbitrary positive constant. Note that the functions $\eta_{1,t}$ and $\eta_{2,t}$ correspond to $\eta_1(t)$ and $\eta_2(t)$, respectively. Since, by assumption, $y_1(t) = y(t) = y_2(t)$ and $\eta_1(0) = \eta_0 = \eta_2(0)$ and by the uniqueness of the mild solution of (\ref{NewDynamicsEta}), the relation $\eta_{1,t} = \eta_{2,t}$ holds true. Consequently,
\begin{align*}
&\Vert\eta_1(\tilde{t})-\eta_2(\tilde{t})\Vert \leq \int_t^{\tilde{t}} \Vert T_Q(\tilde{t}-s)R(y_1(s)-y_2(s))\Vert ds\\
&\hspace{0cm}+\int_t^{\tilde{t}}\Vert T_Q(\tilde{t}-s)P^\perp(\tilde{f}(y_1(s),\eta_1(s))-\tilde{f}(y_2(s),\eta_2(s)))\Vert ds.
\end{align*}
Assumption \ref{OpF} together with the boundedness of the operator $R$ implies that
\begin{align*}
&\Vert\eta_1(\tilde{t})-\eta_2(\tilde{t})\Vert\\
&\leq(\Vert R\Vert_{\mathcal{L}(\mathbb{R},\mathcal{I})}+2l_1)\int_t^{\tilde{t}}Me^{\vert\omega\vert(\tilde{t}-s)} \vert y_1(s)-y_2(s)\vert ds\\
&\hspace{2.5cm}+2l_2\int_t^{\tilde{t}} Me^{\vert\omega\vert(\tilde{t}-s)}\Vert\eta_1(s)-\eta_2(s)\Vert ds,
\end{align*}
where $l_1$ and $l_2$ are the positive constants introduced in (\ref{LipschitzConstants}). We shall use the notation $\Vert R\Vert_{\mathcal{L}(\mathbb{R},\mathcal{I})}+2l_1 =: \mathfrak{g}$ in what follows. Applying Gronwall's lemma to the function $e^{-\vert\omega\vert\tilde{t}}\Vert\eta_1(\tilde{t})-\eta_2(\tilde{t})\Vert$ yields the inequality
\begin{align*}
&\Vert\eta_1(\tilde{t})-\eta_2(\tilde{t})\Vert\\
&\hspace{1.5cm}\leq \mathfrak{g}Me^{\vert\omega\vert\tilde{t}}e^{2Ml_2(\tilde{t}-t)}\int_t^{\tilde{t}}e^{-\vert\omega\vert s}\vert y_1(s)-y_2(s)\vert ds.
\end{align*}
Taking the supremum over all $\tilde{t}$ in $[t,t+\tau]$ on both sides yields the estimate
\begin{align}
&\sup_{\tilde{t}\in[t,t+\tau]}\Vert\eta_1(\tilde{t})-\eta_2(\tilde{t})\Vert\nonumber\\
&\leq \mathfrak{g}Me^{(\vert\omega\vert+2Ml_2)\tau}\tau\sup_{\tilde{t}\in[t,t+\tau]}\vert y_1(\tilde{t})-y_2(\tilde{t})\vert.\label{EstimateLipschitzEta}
\end{align}
The notation $\mathfrak{f}_{\tau} := \mathfrak{g}Me^{(\vert\omega\vert+2Ml_2)\tau}\tau$ will be used for the sake of simplicity. According to the definition (\ref{EquationT}) of the nonlinear operator $T$, it holds that
\begin{align*}
&\vert T(y_1)(\tilde{t})-T(y_2)(\tilde{t})\vert \leq \vert P_0\vert\vert y_1(\tilde{t})-y_2(\tilde{t})\vert\\
&+ \Vert S\Vert_{\mathcal{L}(\mathcal{I},\mathbb{R})}\Vert\eta_1(\tilde{t})-\eta_2(\tilde{t})\Vert + l_1\Vert c\Vert_H\vert y_1(\tilde{t})-y_2(\tilde{t})\vert\\
&+ l_2\Vert c\Vert_H\Vert\eta_1(\tilde{t})-\eta_2(\tilde{t})\Vert,
\end{align*}
which, combined with (\ref{EstimateLipschitzEta}), leads to
\begin{equation*}
\sup_{\tilde{t}\in[t,t+\tau]}\vert T(y_1)(\tilde{t})-T(y_2)(\tilde{t})\vert\leq \rho\sup_{\tilde{t}\in[t,t+\tau]}\vert y_1(\tilde{t})-y_2(\tilde{t})\vert,
\end{equation*}
where $\rho := \vert P_0\vert + \Vert S\Vert_{\mathcal{L}(\mathcal{I},\mathbb{R})}\mathfrak{f}_\tau + l_1\Vert c\Vert_H + l_2\Vert c\Vert_H \mathfrak{f}_\tau$.
\qed
\end{pf}
This means that funnel control is feasible for a nonlinear infinite-dimensional system of the form (\ref{AbstractODE}) which satisfies Assumptions \ref{OpA}, \ref{OpF}, \ref{ShapeFct} and \ref{BIBO_Assum}. Moreover the closed-loop system which consists of the interconnection of (\ref{AbstractODE}), described by (\ref{NewDynamicsY})--(\ref{NewDynamicsEta}) (system $\tilde{\Sigma}$), with the funnel controller (\ref{FunnelController}) has the properties described in Theorem \ref{ThmFunnel}. This system is depicted in Figure \ref{fig:closed-loop}.

We state hereafter a useful criterion for checking Assumption \ref{BIBO_Assum}, assuming that the nonlinearity is bounded.

\begin{prop}\label{BIBO_Stability_Sol}
Assuming that the semigroup $(T_Q(t))_{t\geq 0}$ is exponentially stable on $\mathcal{I}$ and that the nonlinear operator $f$ satisfies $\Vert f(x)\Vert_H\leq \hat{\sigma}$, for some constant $\hat{\sigma} > 0$ independent of $x$, is sufficient to ensure that Assumption \ref{BIBO_Assum} is satisfied.
\end{prop}

\begin{pf}
Let us fix $k_1>0, k_2>0, y\in\mathcal{C}(\mathbb{R}^+,\mathbb{R})$ and $\eta_0\in\mathcal{I}$ such that $\sup_{t\in\mathbb{R}^+}\vert y(t)\vert\leq k_1$ and $\Vert\eta_0\Vert\leq k_2$. As the semigroup $(T_Q(t))_{t\geq 0}$ is exponentially stable, the inequality $\Vert T_Q(t)\Vert\leq \tilde{M}e^{-\tilde{\omega}t}$ holds for some $\tilde{M}\geq 1$ and $\tilde{\omega}>0$. Moreover the mild solution of (\ref{NewDynamicsEta}) with initial condition $\eta_0$ and with the function $y$ fixed above is given by
\begin{align*}
\eta(t) &= T_Q(t)\eta_0 + \int_0^t T_Q(t-s)Ry(s)ds\\
&\hspace{3cm} + \int_0^tT_Q(t-s)P^\perp\tilde{f}(y(s),\eta(s))ds.
\end{align*}
By taking the $H-$norm (restricted to $\mathcal{I}$) on both sides and by using the exponential stability of $(T_Q(t))_{t\geq 0}$, one gets that
\begin{align*}
\Vert\eta(t)\Vert&\leq \tilde{M}e^{-\tilde{\omega} t}\Vert\eta_0\Vert + \int_0^t \tilde{M}e^{-\tilde{\omega}(t-s)}\Vert Ry(s)\Vert ds\\
&\hspace{2cm}+\int_0^t\tilde{M}e^{-\tilde{\omega}(t-s)}\Vert P^\perp\tilde{f}(y(s),\eta(s))\Vert ds.
\end{align*}
The boundedness of the operator $R$ combined with the definition of $P^\perp$ and the assumption on the operator $f$ entails that
\begin{align*}
\Vert\eta(t)\Vert&\leq \tilde{M}e^{-\tilde{\omega} t}\Vert\eta_0\Vert + \tilde{M}e^{-\tilde{\omega} t}\Vert R\Vert_{\mathcal{L}(\mathbb{R},\mathcal{I})}\int_0^t e^{\tilde{\omega} s}\vert y(s)\vert ds\\
&\hspace{4cm}+2\hat{\sigma} \tilde{M}e^{-\tilde{\omega} t}\int_0^t e^{\tilde{\omega} s}ds.
\end{align*}
It follows, by using the estimate $\sup_{t\in\mathbb{R}^+}\vert y(t)\vert\leq k_1$, that
\begin{align*}
\Vert\eta(t)\Vert\leq \tilde{M}e^{-\tilde{\omega} t}\left(\Vert\eta_0\Vert + (\Vert R\Vert_{\mathcal{L}(\mathbb{R},\mathcal{I})}k_1+2\hat{\sigma})\int_0^t e^{\tilde{\omega} s}ds\right).
\end{align*}
Hence
\begin{align}
\Vert\eta(t)\Vert&\leq \tilde{M}\left(e^{-\tilde{\omega} t}\Vert\eta_0\Vert + \frac{\Vert R\Vert_{\mathcal{L}(\mathbb{R},\mathcal{I})}k_1+2\hat{\sigma}}{\tilde{\omega}}(1-e^{-\tilde{\omega} t})\right)\nonumber\\
&\leq \kappa,
\label{EstimEta}
\end{align}
where $\kappa := \tilde{M}\left(k_2 + \frac{\Vert R\Vert_{\mathcal{L}(\mathbb{R},\mathcal{I})}k_1+2\hat{\sigma}}{\tilde{\omega}}\right)$ does not depend on $t$.
\qed
\end{pf}

We conclude this section with an analysis in line with Remark \ref{RemarkShapeFct}, in order to assess the impact of the approximation of the vectors $b$ and $c$ by elements of $D(A)$ and $D(A^*)$, respectively. We shall study here a bound on the error resulting from the approximation. Let us first remember that funnel control aims to be applied on a nonlinear infinite-dimensional system of the form (\ref{AbstractODE}). According to Remark \ref{RemarkShapeFct}, if $b$ and $c$ are not in $D(A)$ and $D(A^*)$, respectively, one may select two sequences $\{b_n\}_{n\in\mathbb{N}}\subset D(A)$ and $\{c_n\}_{n\in\mathbb{N}}\subset D(A^*)$ which converge towards $b$ and $c$ in $H-$norm, respectively. Let us select one element of each sequence, denoted by $b_\ell$ and $c_\ell, \ell\in\mathbb{N}$. For these elements, an approximation of (\ref{AbstractODE}) without any disturbance takes the form
\begin{equation}
\left\{
\begin{array}{l}
\dot{x}_\ell(t) = Ax_\ell(t) + f(x_\ell(t)) + b_\ell u_\ell(t),\\
y_\ell(t) = \langle x_\ell(t),c_\ell\rangle_H,\\
x_\ell(0) = x_0\in H,
\end{array}\right.
\label{AbstractODEApprox}
\end{equation}
where the initial condition $x_\ell(0)$ is purposely chosen the same as the one of (\ref{AbstractODE}). According to the arguments that are developed in this section, funnel control is feasible for (\ref{AbstractODEApprox}). A funnel controller for (\ref{AbstractODEApprox}) is written as
\begin{equation}
u_\ell(t) = \frac{-e_\ell(t)}{1-\phi^2(t)e_\ell^2(t)}, e_\ell(t) = y_\ell(t) - y_{ref}(t).
\label{FunnelControlApprox}
\end{equation}
Moreover, according to Theorem \ref{ThmFunnel}, the control produced by controller (\ref{FunnelControlApprox}) is bounded on $\mathbb{R}^+$, i.e. there exists $k_\ell>0$ such that $\vert u_\ell(t)\vert\leq k_\ell$ for all $t\in\mathbb{R}^+$. Applying this control action to (\ref{AbstractODE}) without taking disturbances into account results in the closed-loop system
\begin{equation*}
\left\{
\begin{array}{l}
\dot{x}(t) = Ax(t) + f(x(t)) + bu_\ell(t),\\
y(t) = \langle x(t),c\rangle_H,\\
x(0) = x_0\in H.
\end{array}\right.
\end{equation*}
Let us study the evolution of the state and output approximation errors $\rho(t) := x(t) - x_\ell(t)$ and $\lambda(t) := y(t) - y_\ell(t)$. The time derivative of the variable $\rho(t)$ is subject to the following abstract ODE
\begin{equation*}
\dot{\rho}(t) = A\rho(t) + f(x(t)) - f(x_\ell(t)) + (b - b_\ell)u_\ell(t), \rho(0) = 0,
\end{equation*}
which has a mild solution given by
\begin{align}
\rho(t) &= \int_0^t T(t-s)[f(x(s))-f(x_\ell(s))]ds\nonumber\\
&+ \int_0^t T(t-s)[b - b_\ell]u_\ell(s)ds,\label{ErrorApprox}
\end{align}
for any positive $t$, where $(T(t))_{t\geq 0}$ denotes the $C_0-$semigroup whose operator $A$ is the infinitesimal generator. Taking the $H-$norm of both sides of (\ref{ErrorApprox}) yields the inequality
\begin{align*}
\Vert\rho(t)\Vert_H &\leq \int_0^t M^*e^{\omega^* (t-s)}l_f\Vert\rho(s)\Vert_H ds\\
&+ \int_0^t M^*e^{\omega^* (t-s)}\Vert b-b_\ell\Vert_H\vert u_\ell(s)\vert ds,
\end{align*}
where the inequality $\Vert T(t)\Vert\leq M^*e^{\omega^* t}, M^*\geq 1, \omega^* >\omega_0$ has been used, $\omega_0$ denotes the growth bound of $(T(t))_{t\geq 0}$, and the positive constant $l_f$ is a Lipschitz constant associated to $f$. By using the boundedness of the input $u_\ell(t)$, it follows that
\begin{align*}
\Vert e^{-\omega^* t}\rho(t)\Vert_H &\leq M^*l_f\int_0^t \Vert e^{-\omega^* s}\rho(s)\Vert_H ds\\
&+\Vert b-b_\ell\Vert_H M^*k_\ell \frac{1-e^{-\omega^* t}}{\omega^*}.
\end{align*}
By using Gronwall's inequality and by chosing $\omega^*$ sufficiently large such that $\omega^*>\max\{0,\omega_0\}$, one gets that
\begin{align}
\Vert\rho(t)\Vert_H\leq \frac{M^*k_\ell (e^{\omega^* t}-1)}{\omega^*}e^{M^*l_ft}\Vert b - b_\ell\Vert_H, t\geq 0.\label{InequalityErrorGronwall}
\end{align}
In view of the previous inequality, for any fixed time instant $t$, letting $\ell$ tend to $\infty$ implies that $x_\ell(t)$ converges towards $x(t)$ in $H-$norm, since $\lim_{\ell\to\infty}\Vert b - b_\ell\Vert_H = 0$, i.e. one gets the pointwise convergence of the sequence of functions $\{x_\ell\}_{\ell\in\mathbb{N}}$ towards the state trajectory $x$. Actually it follows from inequality (\ref{InequalityErrorGronwall}) that this convergence is uniform on any compact interval in $\mathbb{R}^+$. Let us have a look at the approximated output $y_\ell(t)$. For this purpose, let us consider the error function $\lambda(t) := y(t) - y_\ell(t), t\geq 0$. Observe that
\begin{align*}
\vert \lambda(t)\vert &= \vert \langle x(t),c\rangle - \langle x_\ell(t),c_\ell\rangle\vert\\
&=\vert \langle x(t),c\rangle - \langle x_\ell(t),c\rangle + \langle x_\ell(t),c\rangle - \langle x_\ell(t),c_\ell\rangle\vert\\
&=\vert\langle \rho(t),c\rangle + \langle x_\ell(t),c-c_\ell\rangle\vert\\
&\leq \Vert \rho(t)\Vert_H\Vert c\Vert_H + \Vert x_\ell(t)\Vert_H\Vert c - c_\ell\Vert_H,
\end{align*}
where the Cauchy-Schwarz inequality has been used. Fixing the time and letting $\ell$ go to $\infty$ implies that the output error made by the approximation of $b$ and $c$, $\lambda(t)$, goes to\footnote{The relations $\lim_{\ell\to\infty}\Vert c - c_\ell\Vert_H = 0$ and $\lim_{\ell\to\infty}\Vert x(t) - x_\ell(t)\Vert_H = 0$ have been used.} $0$ pointwise and uniformly on compact intervals.

\section{Applications}\label{Simus}
Two applications are considered in this section. First we first shall use the framework developed in the previous section to regulate the average temperature in a nonisothermal plug-flow tubular reactor. As a second application, a nonlinear damped sine-Gordon equation is considered.

\subsection{Nonisothermal plug-flow tubular reactor}
The dimensionless nonlinear dynamics of a plug-flow tubular reactor without axial dispersion is governed by the following PDE
\begin{equation}
\left\{
\begin{array}{l}
\partial_t\theta_1 = -\partial_z\theta_1 + \delta f(\theta_1,\theta_2) + \beta(1_{[0,1]}(z)\theta_w(t)-\theta_1)\\
\partial_t\theta_2 = -\partial_z\theta_2 + f(\theta_1,\theta_2)\\
\theta(t,0) = 0, \theta_2(t,0) = 0,
\end{array}\right.
\label{DynamicsReactor}
\end{equation}
with $t\geq 0, z\in [0,1]$ and where $\theta_1$ and $\theta_2$ represent the dimensionless temperature and concentration inside the reactor, respectively. The nonlinear part of the equation, encompassed in the function $f$, is due to the reaction kinetics. The Arrhenius law is considered here, which implies that that the function $f$ is expressed as
\begin{equation}
f(\theta_1,\theta_2) = \left\{\begin{array}{l}
0 \text{ if } \theta_1 < -1\\
\alpha e^{\frac{\mu\theta_1}{1+\theta_1}} \text{ if } \theta_1\geq -1 \text{ and } \theta_2 < 0\\
\alpha(1-\theta_2)e^{\frac{\mu\theta_1}{1+\theta_1}} \text{ if } \theta_1\geq -1 \text{ and } 0\leq\theta_2\leq 1\\
0 \text{ otherwise. }
\end{array}\right.
\label{NonlinearF}
\end{equation}
Note that this definition of $f$ implies that the latter is uniformly Lipschitz as a pointwise function defined from $\mathbb{R}^2$ into $\mathbb{R}$. Moreover, it satisfies $\vert f(x,y)\vert\leq \alpha e^\mu $ for any $(\begin{matrix} x & y\end{matrix})^T\in\mathbb{R}^2$. The positive constants $\alpha, \beta, \mu$ and $\delta$ depend on the model parameters. An overview of these constants and their physical meaning can be found in \cite{DochainChemical} among others. The scalar control variable, denoted by $\theta_w(t)$, is due to a heat exchanger that acts as a distributed control input along the reactor through the characteristic function $1_{[0,1]}(z)$ which takes the value $1$ for $z\in[0,1]$ and $0$ elsewhere. The control objective that is pursued here consists in the tracking of the following output function which corresponds to the mean value of the temperature along the reactor:
\begin{equation}
y(t) = \int_0^1\theta_1(t,z)dz.
\label{OutputReactor}
\end{equation}
In order to reach this goal, we shall develop a funnel controller producing an input $\theta_w(t)$ which will take the form (\ref{FunnelController}) for some reference signal $y_{\text{ref}}\in W^{1,\infty}(\mathbb{R}^+,\mathbb{R})$. First observe that (\ref{DynamicsReactor}) admits the abstract representation 
\begin{equation}
\left\{
\begin{array}{l}
\dot{x}(t) = Ax(t) + F(x(t)) + Bu(t)\\
x(0) = x_0\in X\\
y(t) = \langle c,x(t)\rangle_X,
\end{array}\right.
\label{AbstractODEExample}
\end{equation}
where the state space $X$ is chosen as $X := L^2(0,1)\times L^2(0,1)$ equipped with the inner product
\begin{equation}
\left\langle\left(\begin{smallmatrix}
x_1\\ 
x_2
\end{smallmatrix}\right),\left(\begin{smallmatrix}
w_1\\ 
w_2
\end{smallmatrix}\right)\right\rangle_X := \langle x_1,w_1\rangle_{L^2(0,1)} + \langle x_2,w_2\rangle_{L^2(0,1)},
\label{InnerProductL2xL2}
\end{equation}
for $(\begin{matrix}
x_1 & x_2
\end{matrix})^T, (\begin{matrix}
w_1 & w_2
\end{matrix})^T\in X$. Note that the inner product on $L^2(0,1)$ (as a real Hilbert space) is defined by
\begin{equation}
\langle x_1,w_1\rangle_{L^2(0,1)} = \int_0^1x_1(z)w_1(z)dz,
\label{InnerProductL2Real}
\end{equation}
where the functions are supposed to take values in $\mathbb{R}$. The unbounded linear operator $A$ is given by $A = \left(\begin{smallmatrix}-d_z - \beta I & 0\\
0 & -d_z\end{smallmatrix}\right)$ on the dense linear subspace
\begin{align*}
D(A)&=\left\{(\begin{matrix}x_1 & x_2\end{matrix})^T\in H^1(0,1)\times H^1(0,1),\right.\\
&\hspace{3cm}\left. x_1(0) = 0 = x_2(0)\right\}.
\end{align*}
From that definition of $(A,D(A))$, the operator $A$ is the infinitesimal generator of a contraction (and even exponentially stable) $C_0-$semigroup on $X$. Hence Assumption \ref{OpA} is satisfied. The nonlinear operator $F : X\to X$ is given by $F(x_1,x_2) = (\begin{matrix}\delta f(x_1,x_2) & f(x_1,x_2)\end{matrix})^T$. Due to the definition of $f$ given in (\ref{NonlinearF}) and the properties of that function viewed as a pointwise function acting on $\mathbb{R}^2$, the nonlinear operator $F$ is uniformly Lipschitz continuous on $X$. Moreover, it satisfies $\Vert F(x_1,x_2)\Vert_X\leq \alpha e^\mu \sqrt{\delta^2+1}$ for any $(\begin{matrix}x_1 & x_2\end{matrix})^T\in X$. This implies that Assumption \ref{OpF} is fulfilled. The control operator $B:\mathbb{R}\to X$ is given by $Bu = (\begin{matrix}\beta 1_{[0,1]}(z) & 0\end{matrix})^Tu$ while the observation operator $C:X\to\mathbb{R}$ takes the form
\begin{equation*}
C(\begin{matrix} x_1 & x_2\end{matrix})^T = \left\langle \left(\begin{matrix}1_{[0,1]}\\ 0\end{matrix}\right),\left(\begin{matrix}x_1\\ x_2\end{matrix}\right)\right\rangle_X.
\end{equation*}
Hence the functions $b$ and $c$ satisfy $b(z) = \beta c(z) = (\begin{matrix}\beta 1_{[0,1]}(z) & 0\end{matrix})^T$. Note that the domain $D(A^*)$ of the adjoint operator of $A$ is given by 
\begin{equation*}
\left\{(\begin{matrix}x_1 & x_2\end{matrix})^T\in H^1(0,1)\times H^1(0,1), x_1(1) = 0 = x_2(1)\right\}.
\end{equation*}
Clearly $b$ and $c$ do not lie in $D(A)$ and $D(A^*)$, respectively. In order to overcome this difficulty, we approximate the function $1_{[0,1]}(z)$ by a linear combination of elements of the orthonormal basis $\left\{\sqrt{2}\sin(n\pi z)\right\}_{n\in\mathbb{N}_0}$ of $L^2(0,1)$. For a fixed $N\in\mathbb{N}_0$, the $N-$th order approximation of $1_{[0,1]}(z)$, denoted by $1_N(z)$, is given by
\begin{align}
1_N(z) &= \sum_{n=1}^N \langle 1_{[0,1]}(\cdot),\sqrt{2}\sin(n\pi\cdot)\rangle_{L^2(0,1)}\sqrt{2}\sin(n\pi z)\nonumber\\
&= \sum_{\overset{n=1}{n \text{ odd}}}^N \frac{4}{n\pi}\sin(n\pi z).\label{1_N}
\end{align}
As this approximation lies in $H^1(0,1)$ and vanishes both for $z = 0$ and $z=1$, the approximations of $b(z)$ and $c(z)$, denoted by $b_N(z)$ and $c_N(z)$ and whose expressions are given by $(\begin{matrix}\beta 1_N(z) & 0\end{matrix})^T$ and $(\begin{matrix}1_N(z) & 0\end{matrix})^T$ are in $D(A)$ and $D(A^*)$, respectively. Now observe that the inner product between $b_N$ and $c_N$ is given by
\begin{align*}
\langle b_N,c_N\rangle_X = \beta \langle 1_N, 1_N\rangle_{L^2(0,1)} = \beta\sum_{\overset{n=1}{n \text{ odd}}}^N \frac{8}{n^2\pi^2} > 0,
\end{align*}
which implies that Assumption \ref{ShapeFct} is satisfied. Now it remains to show that the system $\Sigma_{y\to\eta}$ is BIBO stable in the sense of Assumption \ref{BIBO_Assum}. We shall use the caracterization of Proposition \ref{BIBO_Stability_Sol} to reach this aim. To this end, let us first introduce the decomposition of the state space $X$ as in (\ref{DecompositionH}), which is given here by:
\begin{align*}
X &= \text{span}\left\{c_N\right\}\oplus \left\{b_N\right\}^\perp = \text{span}\left\{c_N\right\}\oplus \left\{c_N\right\}^\perp \\
&= \text{span}\left\{\left(\begin{smallmatrix}1_N\\0\end{smallmatrix}\right)\right\}\oplus \left\{\left(\begin{smallmatrix}
x_1\\x_2
\end{smallmatrix}\right)\in X, \left\langle \left(\begin{smallmatrix}
x_1\\x_2
\end{smallmatrix}\right), \left(\begin{smallmatrix}
1_N\\0
\end{smallmatrix}\right)\right\rangle_X = 0\right\}\\
&= \mathcal{O}\oplus\mathcal{I},
\end{align*}
where $\mathcal{I}$ may also be written as
\begin{equation}
\mathcal{I} = \left\{x_1\in L^2(0,1), \langle x_1,1_N\rangle_{L^2(0,1)} = 0\right\}\times L^2(0,1).
\label{Inputsubspace}
\end{equation}
The operator $Q: D(Q) = D(A)\cap \mathcal{I}\to \mathcal{I}$ whose definition is given in (\ref{OpQ}) is expressed as
\begin{equation*}
Q\eta = P^\perp A\eta - \frac{\langle c_N,\eta\rangle_X}{\langle c_N,b_N\rangle_X}P^\perp Ab_N,
\end{equation*}
for $\eta = (\begin{matrix}\eta_1 & \eta_2\end{matrix})^T\in D(Q)$. Expanding the application of $Q$ to $\eta$ gives rise to
\begin{align*}
Q\eta &= A\eta - \frac{\langle A\eta,b_N\rangle_X}{\langle b_N,b_N\rangle_X}b_N\\
&= \left(\begin{matrix}
-d_z\eta_1 - \beta\eta_1\\
-d_z\eta_2
\end{matrix}\right) + \frac{\langle d_z\eta_1,1_N\rangle_{L^2(0,1)}}{\langle 1_N,1_N\rangle_{L^2(0,1)}}\left(\begin{matrix}
1_N\\
0
\end{matrix}\right)\\
&=: \left(\begin{matrix}
Q_1 & 0\\
0 & Q_2
\end{matrix}\right)\left(\begin{matrix}
\eta_1\\
\eta_2
\end{matrix}\right),
\end{align*}
where the relation $\langle \eta_1,1_N\rangle_{L^2(0,1)} = 0$ has been used. The operator $Q_1: D(Q_1)\subset \left\{1_N\right\}^\perp\to\left\{1_N\right\}^\perp$ is defined as
\begin{align*}
Q_1\eta_1 = -d_z\eta_1 - \beta\eta_1 + \frac{\langle d_z\eta_1,1_N\rangle_{L^2(0,1)}}{\langle 1_N,1_N\rangle_{L^2(0,1)}}1_N,
\end{align*}
for $\eta_1\in D(Q_1)$ given by
\begin{align*}
D(Q_1) = \left\{x_1\in H^1(0,1), x_1(0) = 0\right\}\cap\left\{1_N\right\}^\perp,
\end{align*}
whereas $Q_2\eta_2 = -d_z\eta_2$ on 
\begin{align*}
D(Q_2) = \left\{x_2\in H^1(0,1), x_2(0) = 0\right\}.
\end{align*}
According to Lyapunov's Theorem, see e.g. \cite[Theorem 4.1.3]{CurtainZwartNew}, the semigroup generated by the operator $Q$ is exponentially stable if and only if there exists a positive self-adjoint operator $P\in\mathcal{L}(\mathcal{I})$ such that
\begin{equation}
\langle Q\eta,P\eta\rangle_\mathcal{I} + \langle P\eta,Q\eta\rangle_\mathcal{I} \leq -\langle\eta,\eta\rangle_\mathcal{I},
\label{Lyapunov}
\end{equation}
for all $\eta\in D(Q)$. Note that, for any two functions in $\mathcal{I}$, their inner product on $\mathcal{I}$ is the same as the one defined in (\ref{InnerProductL2xL2}). Let us propose the following form for the operator $P$
\begin{align*}
P = \left(\begin{matrix}
P_1 & 0\\
0 & P_2
\end{matrix}\right).
\end{align*}
Let us define $P_1$ as $P_1\eta_1 = \frac{1}{2\beta}\eta_1$ and $P_2$ as
\begin{align*}
(P_2\eta_2)(z) = (1-z)\eta_2(z), z\in [0,1].
\end{align*}
Now observe that, for $\eta_1\in D(Q_1)$, there holds
\begin{align*}
&\langle Q_1\eta_1,P_1\eta_1\rangle + \langle P_1\eta_1,Q_1\eta_1\rangle = 2\langle Q_1\eta_1,P_1\eta_1\rangle\\
&=-\frac{1}{\beta}\langle d_z\eta_1,\eta_1\rangle - \langle\eta_1,\eta_1\rangle + \frac{1}{\beta}\frac{\langle d_z\eta_1,1_N\rangle_{L^2(0,1)}}{\langle 1_N,1_N\rangle_{L^2(0,1)}}\langle 1_N,\eta_1\rangle\\
&=-\frac{1}{2\beta}\eta_1^2(1) - \langle\eta_1,\eta_1\rangle\leq - \langle\eta_1,\eta_1\rangle,
\end{align*}
where the fact that $\eta_1$ and $1_N$ are orthogonal has been used. Moreover, for any $\eta_2\in D(Q_2)$, one has that
\begin{align*}
&\langle Q_2\eta_2,P_2\eta_2\rangle + \langle P_2\eta_2,Q_2\eta_2\rangle\\
&= -2\int_0^1(1-z)\eta_2 (d_z\eta_2)dz = \left[-(1-z)\eta_2^2(z)\right]_0^1 - \langle\eta_2,\eta_2\rangle\\
&= - \langle\eta_2,\eta_2\rangle.
\end{align*}
Thus the semigroup generated by the operator $Q$ is exponentially stable. This fact combined with the properties of the nonlinear operator $F$ and Proposition \ref{BIBO_Stability_Sol} ensures that Assumption \ref{BIBO_Assum} is satisfied. Hence funnel control is feasible for the system (\ref{DynamicsReactor}), which means that considering a heat exchanger temperature expressed as
\begin{equation}
\theta_w(t) = \frac{-e(t)}{1-\phi^2(t)e^2(t)},
\label{FunnelControllerReactor}
\end{equation}
where $e(t) = y(t) - y_{\text{ref}}(t), y_{\text{ref}}\in W^{1,\infty}(\mathbb{R}^+,\mathbb{R}), \phi\in \Phi$, and $y$ is given by (\ref{OutputReactor}), yields a closed-loop system which possesses the properties stated in Theorem \ref{ThmFunnel}.

As an illustration of the results, we shall report numerical simulations hereafter. The parameters related to the model (\ref{DynamicsReactor}) have been chosen as follows: $\delta = 0.25, \alpha = 2.3248, \mu = 16.6607, \beta = 8$, see \cite{Aksikas2}. The number of basis functions $N$ in the approximation of the function $b$ has been set at $N = 100$. As initial conditions for $\theta_1$ and $\theta_2$ we consider the following functions
\begin{align*}
\theta_1(0,z) &= 0.02(-z^3+z^2+z),\\
\theta_2(0,z) &= 0.7(-z^3+z^2+z).
\end{align*}
The reference signal that has to be tracked by the output (\ref{OutputReactor}) is $y_{\text{ref}}(t) = \frac{1}{20} + \frac{1}{20}\arctan(t)$. The funnel boundary is chosen as $\frac{1}{\phi(t)} = e^{-2t} + 2.5*10^{-3}$.

\begin{figure}[h!]
\begin{center}
\includegraphics[scale=0.45,trim=1.5cm 3cm 0cm 0.1cm]{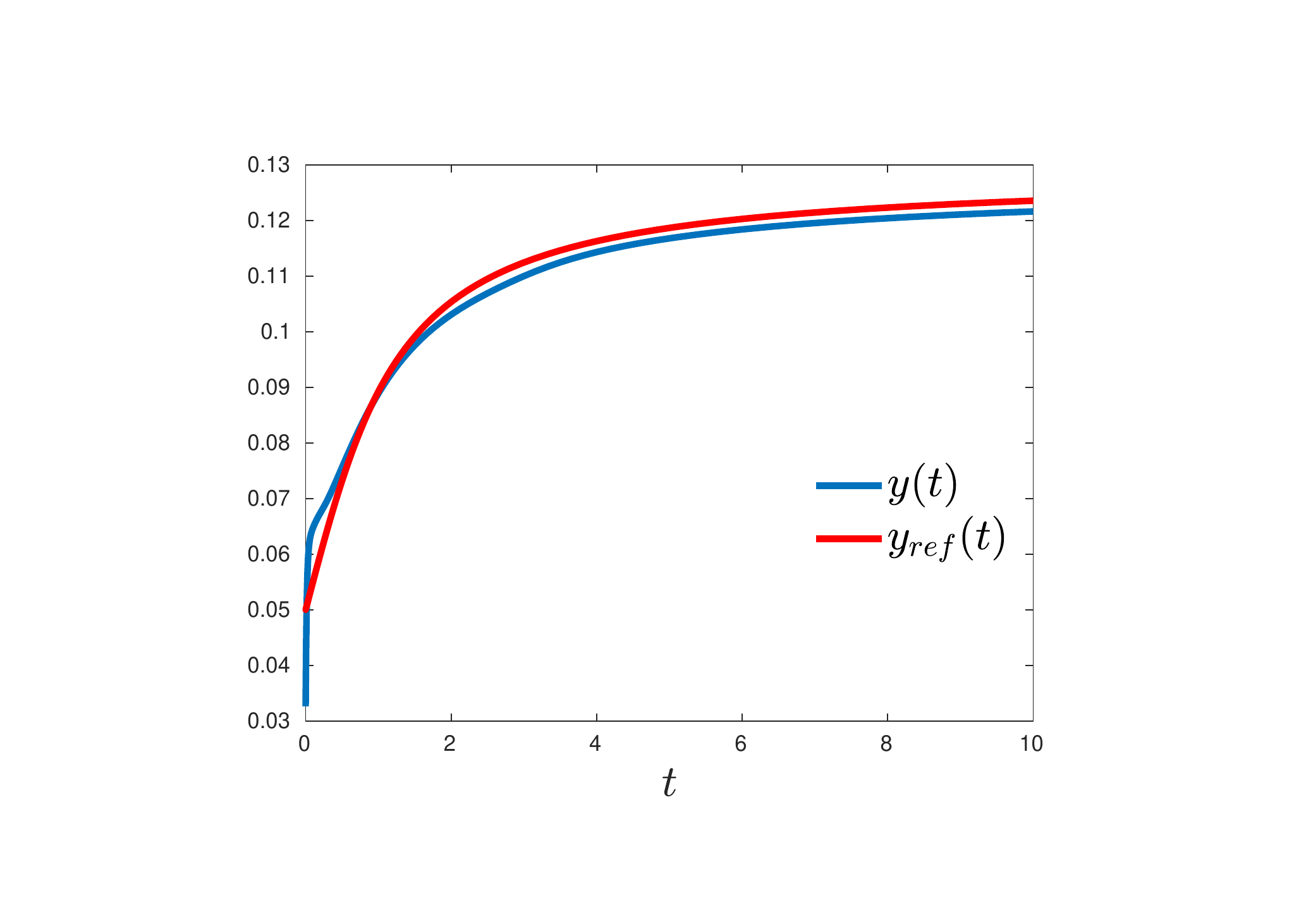}
\end{center}
\caption{Output trajectory (\ref{OutputReactor}) with the reference signal $y_{\text{ref}}(t)$.\label{fig:OutputReactor}}
\end{figure}

\begin{figure}
\begin{center}
\includegraphics[scale=0.45,trim=1.5cm 3cm 0cm 0.1cm]{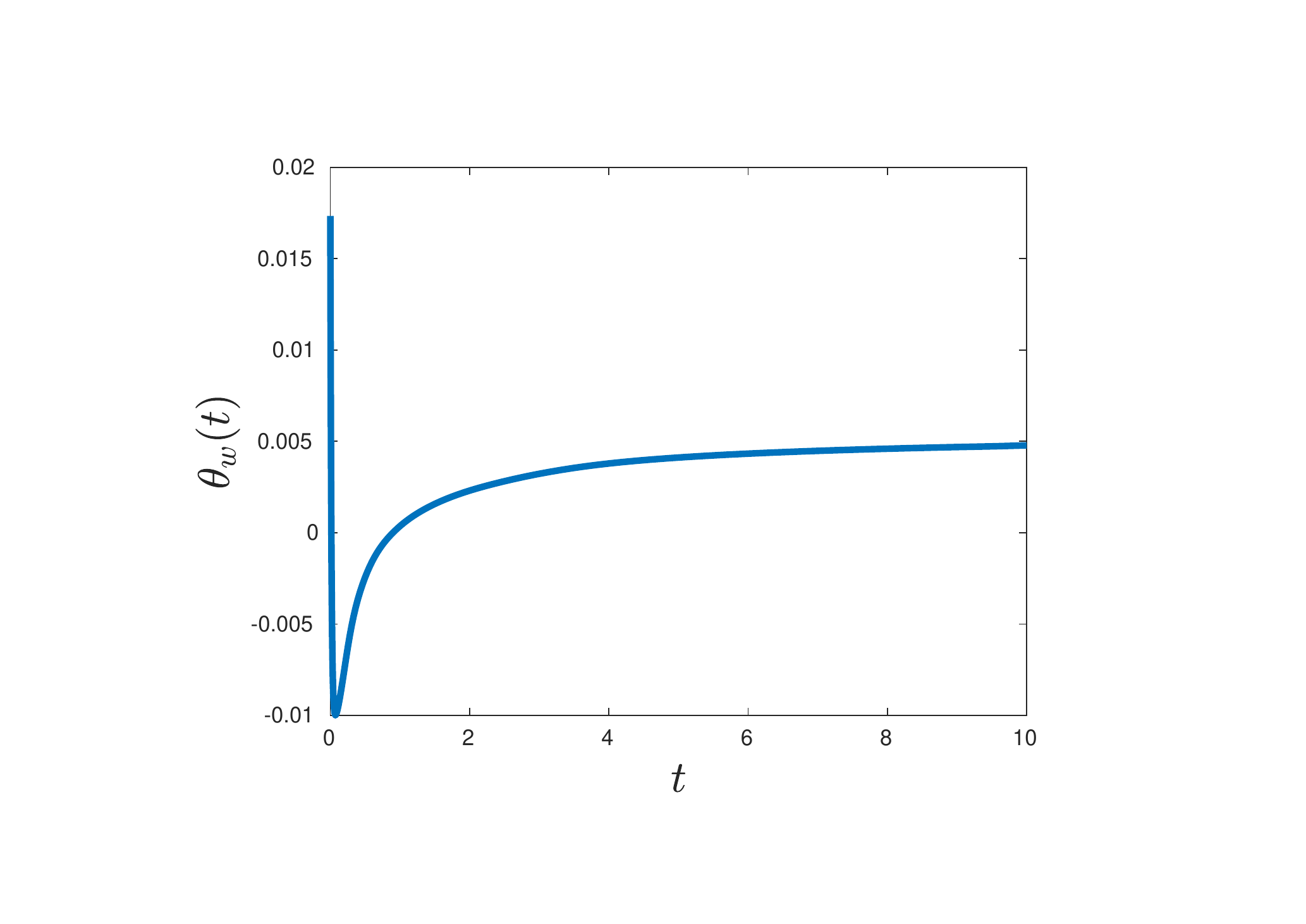}
\end{center}
\caption{Input trajectory.\label{fig:InputReactor}}
\end{figure}

\begin{figure}
\begin{center}
\includegraphics[scale=0.45,trim=1.5cm 3cm 0cm 0.1cm]{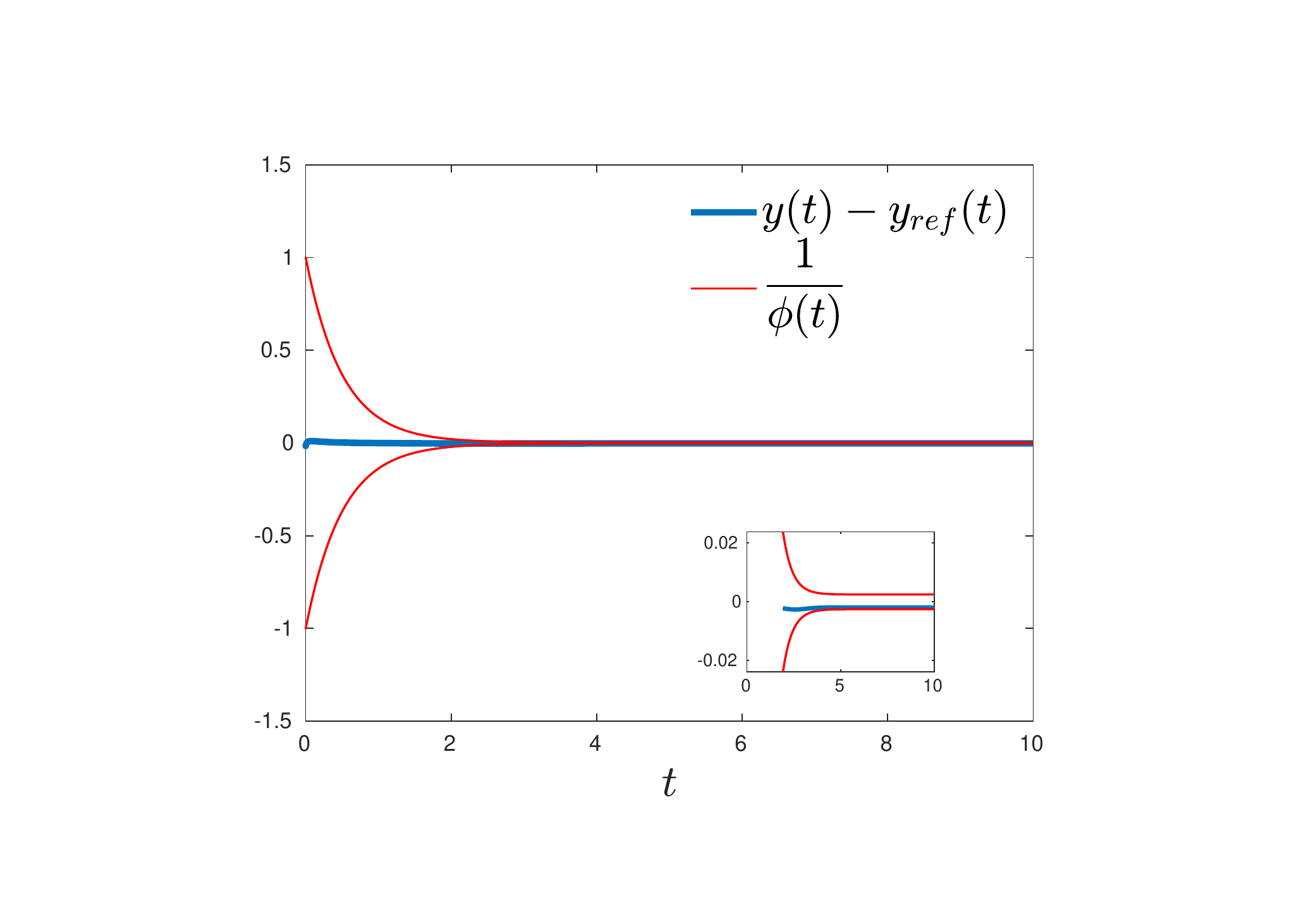}
\end{center}
\caption{Output error trajectory $y(t)-y_{\text{ref}}(t)$ with the funnel boundaries $\frac{1}{\phi(t)}$ and $-\frac{1}{\phi(t)}$.\label{fig:ErrorReactor}}
\end{figure}

\begin{figure}
\begin{center}
\includegraphics[scale=0.45,trim=1.5cm 3cm 0cm 0.1cm]{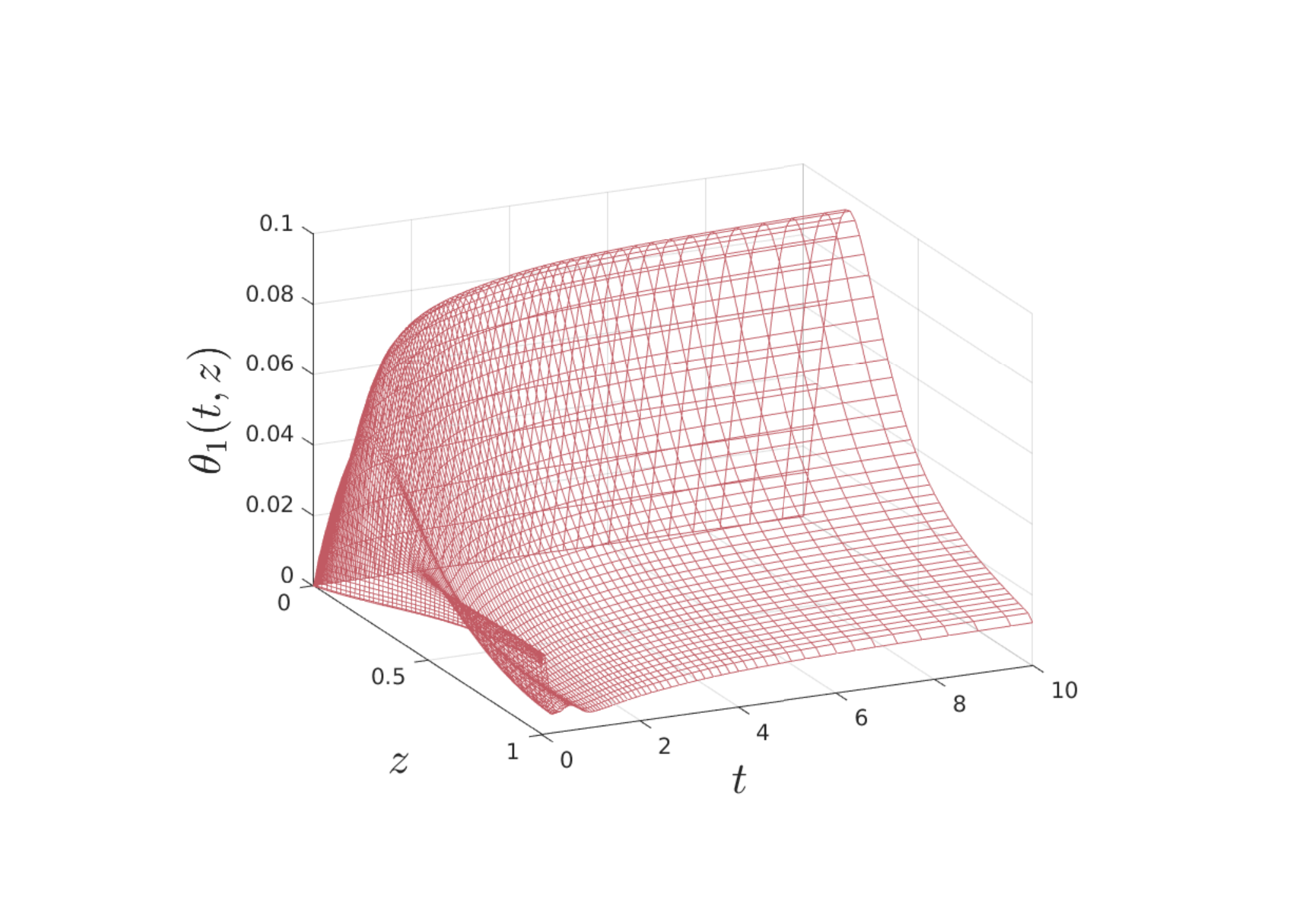}
\end{center}
\caption{Closed-loop state trajectory $\theta_1(t,z)$.\label{fig:tempReactor}}
\end{figure}

\begin{figure}
\begin{center}
\includegraphics[scale=0.45,trim=1.5cm 3cm 0cm 0.1cm]{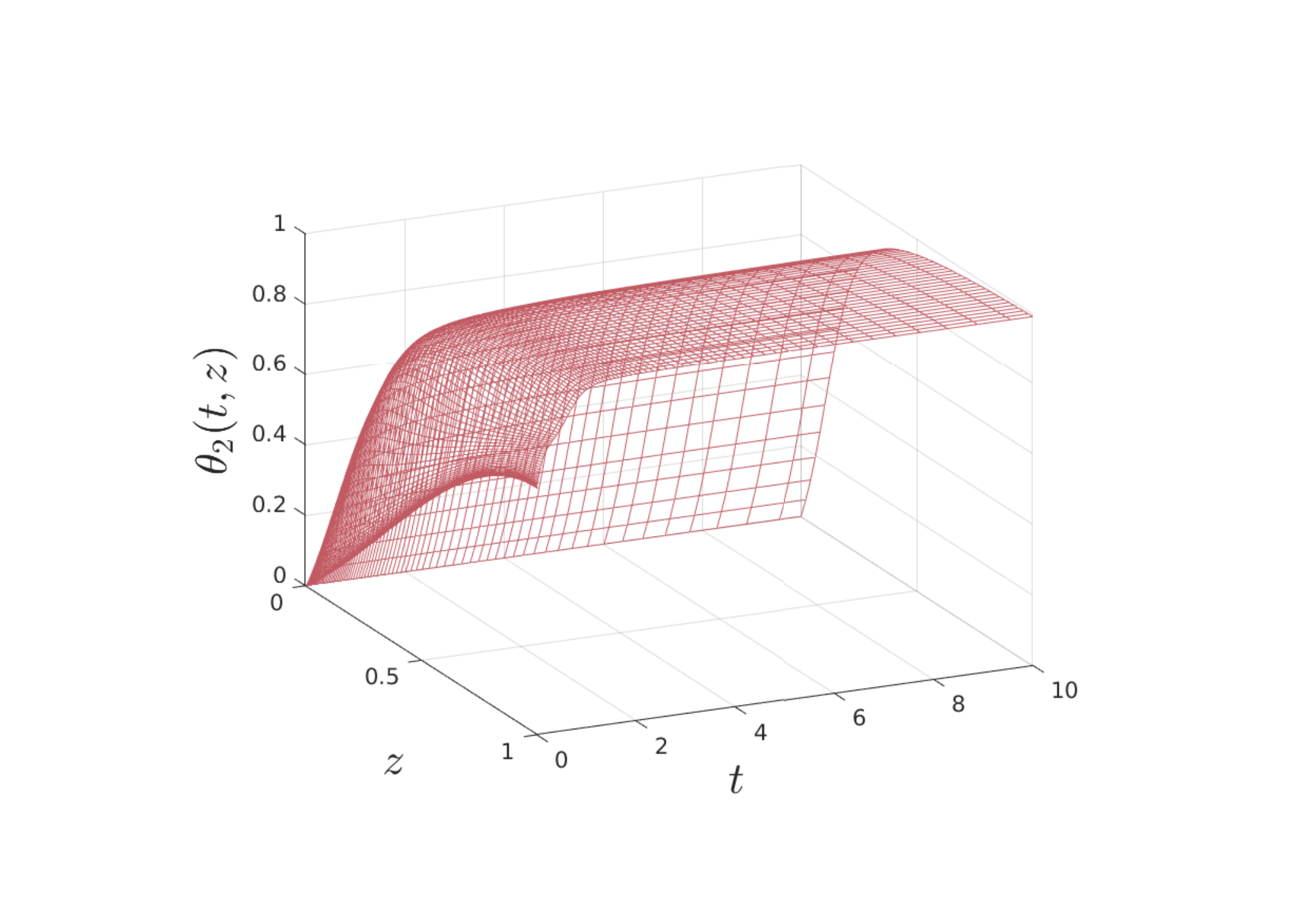}
\end{center}
\caption{Closed-loop state trajectory $\theta_2(t,z)$.\label{fig:concReactor}}
\end{figure}

The output trajectory (\ref{OutputReactor}) and the reference signal $y_{\text{ref}}(t)$ are depicted in Figure \ref{fig:OutputReactor} and the funnel control, $\theta_w(t)$, given by (\ref{FunnelControllerReactor}), is shown in Figure \ref{fig:InputReactor}. The tracking error is depicted in Figure \ref{fig:ErrorReactor}. The state trajectories corresponding to $\theta_1(t,z)$ and $\theta_2(t,z)$ are represented in Figures \ref{fig:tempReactor} and \ref{fig:concReactor}, respectively.

\subsection{Damped sine-Gordon equation}

Here we consider a damped and actuated sine-Gordon equation of the form
\begin{equation}
\left\{
\begin{array}{l}
\partial_{tt}x = \partial_{zz}x - \alpha\partial_t x + \nu \sin(x) + b(z)u(t)\\
x(t,0) = 0, \partial_z x(t,1) = 0,
\end{array}\right.
\label{sineGordon}
\end{equation}
where the space variable $z\in[0,1]$ and $t\in\mathbb{R}^+$ denotes the time variable. The parameters $\nu$ and $\alpha$ are such that $\nu\in\mathbb{R}_0$ and $\alpha> \pi$. This nonlinear PDE encompasses many phenomena in physics as the dynamics of a Josephson junction driven by a current source, see e.g. \cite{Temam,Cuevas}, as well as the dynamics of mechanical transmission lines, see \cite{Cirillo} among others. The stability of the homogeneous dynamics ($u\equiv 0$) of (\ref{sineGordon}) has been investigated for Dirichlet and Neumann boundary conditions in \cite{sineGordonSIAM} and \cite{Callegari}. For control problems related to (\ref{sineGordon}), we refer to \cite{Dolgopolik} and \cite{Efimov} for instance, where boundary energy control and robust input-to-state stability are developed.

Let us consider the operator $A_0 = -\frac{d^2}{dz^2}$ on the domain 
\begin{equation*}
D(A_0) = \left\{x\in H^2(0,1), x(0) = 0 = \frac{dx}{dz}(1)\right\}.
\end{equation*}
As the operator $A_0$ is self-adjoint and coercive\footnote{By using Poincar\'e's inequality it can be seen that, for any $x\in D(A_0)$, the relation $\langle A_0x,x\rangle_{L^2(0,1)}\geq \frac{\pi^4}{4}\Vert x\Vert_{L^2(0,1)}^2$ holds.}, it admits a unique nonnegative square-root, see e.g. \cite[Lemma A.3.82]{CurtainZwartNew}, which satisfies 
\begin{equation*}
D(A_0^\frac{1}{2}) = \left\{x\in H^1(0,1), x(0) = 0\right\}.
\end{equation*}
This allows us to consider the Hilbert state space $Z = D(A_0^{\frac{1}{2}})\times X$ equipped with the inner product\footnote{The equality $\langle d_z\zeta_1,d_z w_1\rangle_X = \langle A_0^{\frac{1}{2}}\zeta_1,A_0^{\frac{1}{2}}w_1\rangle_X$ holds.}
\begin{equation}
\left\langle\left(\begin{matrix}
\zeta_1\\
\zeta_2
\end{matrix}\right),\left(\begin{matrix}
w_1\\
w_2
\end{matrix}\right)\right\rangle_Z := \langle d_z\zeta_1,d_zw_1\rangle_{X} + \langle \zeta_2,w_2\rangle_{X},
\label{InneProductZ}
\end{equation}
where $\zeta_1, w_1\in D(A_0^{\frac{1}{2}})$ and $\zeta_2, w_2\in X$ with $X := L^2(0,1)$. The inner product on $X$ is the same as the one defined in (\ref{InnerProductL2Real}). By considering the state variable $\zeta = \left(\begin{smallmatrix}
x\\\partial_t x\end{smallmatrix}\right) =: \left(\begin{smallmatrix}
\zeta_1\\\zeta_2\end{smallmatrix}\right)$, the PDE (\ref{sineGordon}) may be written as (\ref{AbstractODEExample}), where the operator $A$ is given by
\begin{equation}
A = \left(\begin{matrix}
0 & I\\
-A_0 & -\alpha I
\end{matrix}\right)
\end{equation}
on $D(A) = D(A_0)\times D(A_0^{\frac{1}{2}})$. According to \cite[Example 2.3.5]{CurtainZwartNew}, the operator $A$ is the generator of a contraction $C_0-$semigroup on $Z$. Note that the adjoint operator of $A$, denoted by $A^*$, is expressed as $A^* = \left(\begin{smallmatrix}0 & -I\\
A_0 & -\alpha I\end{smallmatrix}\right)$ on $D(A^*) = D(A)$. The nonlinear operator $F:Z\to Z$ is expressed as $F(\zeta_1,\zeta_2) = \left(\begin{smallmatrix}0\\\nu\sin(\zeta_1)\end{smallmatrix}\right)$. The latter is uniformly Lipschitz continuous and satisfies $\Vert F(\zeta_1,\zeta_2)\Vert_Z\leq \vert\nu\vert$ for any $(\begin{matrix}
\zeta_1 & \zeta_2
\end{matrix})^T\in Z$. Consequently, Assumptions \ref{OpA} and \ref{OpF} are satisfied. Here we consider that the operator $B:\mathbb{R}\to Z$ is defined for $u\in\mathbb{R}$ by $Bu = b(z)u = \left(\begin{smallmatrix}0\\ b_N(z)\end{smallmatrix}\right)u$ where $b_N(z) = 1_N(z) = \sum_{\overset{n=1}{n \text{ odd}}}^N \frac{4}{n\pi}\sin(n\pi z)$ is the function defined in (\ref{1_N}). Obviously $b\in D(A)$. For the function $c$, we choose the expression
\begin{equation*}
c(z) = \left(\begin{matrix}
\frac{2}{\pi^2}(\alpha + \sqrt{\alpha^2-\pi^2})\sin(\frac{\pi}{2}z)\\
\sin(\frac{\pi}{2}z)
\end{matrix}\right),
\end{equation*}
so that the output trajectory corresponding to (\ref{sineGordon}) is given by
\begin{align}
y(t) &= \langle c,\zeta(t)\rangle_Z\nonumber\\
&= \frac{\alpha + \sqrt{\alpha^2-\pi^2}}{\pi}\int_0^1 \cos\left(\frac{\pi}{2}z\right)\partial_z x(t,z)dz\nonumber\\
&\hspace{3cm} + \int_0^1\sin\left(\frac{\pi}{2}z\right)\partial_t x(t,z)dz.
\label{OutputSineGordon}
\end{align}
It can be easily seen that the function $c\in D(A^*)$. A straightforward computation shows that
\begin{equation*}
\langle b,c\rangle_Z = \sum_{\overset{n=1}{n \text{ odd}}}^N \frac{16}{\pi^2(4n^2-1)} > 0,
\end{equation*}
which entails that Assumption \ref{ShapeFct} is satisfied. Before showing that funnel control is feasible in our context, let us introduce the decomposition of the state space 
\begin{align*}
Z = \text{span}\{c\}\oplus \left\{b\right\}^\perp =\mathcal{O}\oplus \mathcal{I},
\end{align*}
with $\mathcal{I} = D(A_0^{\frac{1}{2}})\times \left\{x\in X, \langle x,b_N\rangle_{L^2(0,1)} = 0\right\}$. According to (\ref{ChangeVar}) the system (\ref{sineGordon}) admits the representation (\ref{NewDynamicsY})--(\ref{NewDynamicsEta}), in which we shall focus on the linear part, i.e. the operator $\left(\begin{smallmatrix}P_0 & S\\
R & Q\end{smallmatrix}\right)$. In order to show that funnel control is feasible for (\ref{sineGordon}), one sould use the criterion of BIBO stability stated in Proposition \ref{BIBO_Stability_Sol}. Therefore it remains to show that the $C_0-$semigroup generated by the operator $Q$, see (\ref{OpQ}), is exponentially stable. First observe that the operator $P_0$ takes the form $P_0y = p_0y$, where
\begin{align*}
p_0 &= \frac{\langle A^*c,b\rangle_Z}{\langle b,c\rangle_Z}\\
&= \frac{\langle \frac{1}{2}(\alpha + \sqrt{\alpha^2-\pi^2})\sin\left(\frac{\pi}{2}z\right)-\alpha\sin\left(\frac{\pi}{2}z\right),b_N(z)\rangle_X}{\langle \sin\left(\frac{\pi}{2}z\right),b_N(z)\rangle_X}\\
&= -\frac{\alpha}{2}+\frac{1}{2}\sqrt{\alpha^2-\pi^2} < 0,
\end{align*}
which means that the semigroup generated by $P_0$, which is given by $\left(e^{(-\frac{\alpha}{2}+\frac{1}{2}\sqrt{\alpha^2-\pi^2})t}\right)_{t\geq 0}$ is exponentially stable. Secondly, let us compute the operator $S:\mathcal{I}\to\mathbb{R}, S\eta = \langle\eta,P_\mathcal{I}A^*c\rangle_Z$. Observe that 
\begin{align*}
P_\mathcal{I}A^*c &= A^*c - P_0 c\\
&= \left(\begin{smallmatrix}-1\\
-\frac{\alpha}{2} + \frac{1}{2}\sqrt{\alpha^2-\pi^2}\end{smallmatrix}\right)\sin\left(\frac{\pi}{2}z\right)\\
&\hspace{0cm} + \left(\frac{\alpha}{2}-\frac{1}{2}\sqrt{\alpha^2-\pi^2}\right)\left(\begin{smallmatrix}\frac{2}{\pi^2}(\alpha + \sqrt{\alpha^2-\pi^2})\\
1\end{smallmatrix}\right)\sin\left(\frac{\pi}{2}z\right)\\
&= \left(\begin{smallmatrix}0\\ 0\end{smallmatrix}\right).
\end{align*}
Consequently, the operator $\left(\begin{smallmatrix}P_0 & S\\
R & Q\end{smallmatrix}\right)$ is a triangular operator of the form $\left(\begin{smallmatrix}P_0 & 0\\
R & Q\end{smallmatrix}\right)$. As it is similar to the operator $A$, the corresponding semigroups are also similar, i.e. denoting by $(S(t))_{t\geq 0}$ and by $(\tilde{S}(t))_{t\geq 0}$ the $C_0-$semigroups generated by $A$ and $\left(\begin{smallmatrix}P_0 & S\\
R & Q\end{smallmatrix}\right)$, respectively, the relation $\tilde{S}(t) = U^{-*}S(t)U^*$ holds for all $t\geq 0$. Consequently, $(S(t))_{t\geq 0}$ and $(\tilde{S}(t))_{t\geq 0}$ have the same growth bounds. In that way, let us have a look at the sign of the growth bound of the semigroup $(S(t))_{t\geq 0}$ via a Lyapunov equation approach. Therefore let us define the operator $\Pi\in\mathcal{L}(Z)$ by
\begin{equation*}
\Pi = \left(\begin{matrix}
\frac{1}{\alpha}I + \frac{\alpha}{2}A_0^{-1} & \frac{1}{2}A_0^{-1}\\
\frac{1}{2}I & \frac{1}{\alpha}I
\end{matrix}\right),
\end{equation*}
where $A_0^{-1}$, the inverse of $A_0$, is a bounded and linear operator, since $A_0$ is self-adjoint and coercive, see \cite[Lemma A.3.85]{CurtainZwartNew}. The operator $\Pi$ is self-adjoint for the inner product defined in (\ref{InneProductZ}). Moreover it is coercive since
\begin{equation*}
\langle\Pi \zeta,\zeta\rangle_Z \geq \frac{1}{2\alpha}\Vert \zeta\Vert_Z^2 \text{, for all $\zeta\in Z$. }
\end{equation*}
Furthermore, by taking $\zeta = (\begin{matrix}\zeta_1 & \zeta_2\end{matrix})^T\in D(A)$, one gets that
\begin{align}
&\langle A\zeta,\Pi\zeta\rangle_Z + \langle \Pi\zeta,A\zeta\rangle_Z = 2\langle A\zeta,\Pi\zeta\rangle_Z\nonumber\\
&=2\left\langle \left(\begin{matrix}\zeta_2\\-A_0\zeta_1-\alpha\zeta_2\end{matrix}\right),\left(\begin{matrix}\frac{1}{\alpha}\zeta_1 + \frac{\alpha}{2}A_0^{-1}\zeta_1 + \frac{1}{2}A_0^{-1}\zeta_2\\
\frac{1}{2}\zeta_1 + \frac{1}{\alpha}\zeta_2\end{matrix}\right)\right\rangle_Z\nonumber\\
&= \frac{2}{\alpha}\langle A_0^{\frac{1}{2}}\zeta_2,A_0^{\frac{1}{2}}\zeta_1\rangle_X + \alpha\langle A_0^{\frac{1}{2}}\zeta_2,A_0^{\frac{1}{2}}A_0^{-1}\zeta_1\rangle_X\nonumber\\
& + \langle A_0^{\frac{1}{2}}\zeta_2,A_0^{\frac{1}{2}}A_0^{-1}\zeta_2\rangle_X - \langle A_0\zeta_1,\zeta_1\rangle_X - \frac{2}{\alpha}\langle A_0\zeta_1,\zeta_2\rangle_X\nonumber\\
&-\alpha\langle\zeta_2,\zeta_1\rangle_X - 2\langle\zeta_2,\zeta_2\rangle_X.\label{LyapunovWave}
\end{align}
As $A_0^{-1}\zeta_1$ and $A_0^{-1}\zeta_2$ are in $D(A_0)$ by construction and as the relation $\langle A_0\zeta_1,\zeta_1\rangle_X = \langle A_0^{\frac{1}{2}}\zeta_1,A_0^{\frac{1}{2}}\zeta_1\rangle_X$ holds for any $\zeta_1\in D(A_0)$, the relation (\ref{LyapunovWave}) can also be written as
\begin{align*}
&\langle A\zeta,\Pi\zeta\rangle_Z + \langle \Pi\zeta,A\zeta\rangle_Z = \frac{2}{\alpha}\langle\zeta_2,A_0\zeta_1\rangle_X + \alpha\langle\zeta_2,\zeta_1\rangle_X\\
& + \langle\zeta_2,\zeta_2\rangle_X - \langle A_0^{\frac{1}{2}}\zeta_1,A_0^{\frac{1}{2}}\zeta_1\rangle_X - \frac{2}{\alpha}\langle A_0\zeta_1,\zeta_2\rangle_X\\
&- \alpha\langle\zeta_2,\zeta_1\rangle_X - 2\langle\zeta_2,\zeta_2\rangle_X\\
& = -\langle A_0^{\frac{1}{2}}\zeta_1,A_0^{\frac{1}{2}}\zeta_1\rangle_X - \langle \zeta_2,\zeta_2\rangle_X = -\left\Vert\left(\begin{matrix}\zeta_1\\\zeta_2\end{matrix}\right)\right\Vert_Z^2. 
\end{align*} 
According to \cite[Theorem 4.3.1]{CurtainZwartNew}, it follows that the semigroup $(S(t))_{t\geq 0}$ whose infinitesimal generator is the operator $A$, is exponentially stable. This means that the growth bound of $(S(t))_{t\geq 0}$ is negative. The same holds for the growth bound of the semigroup generated by the operator $\left(\begin{smallmatrix}P_0 & 0\\
R & Q\end{smallmatrix}\right)$. As the growth bound of the semigroup generated by $P_0$ is negative too, the growth bound of the semigroup generated by $Q$ is also negative, showing that this semigroup is exponentially stable. Thanks to Proposition \ref{BIBO_Stability_Sol} the system $\Sigma_{y\to\eta}$ is BIBO stable in the sense of Assumption \ref{BIBO_Assum}. As a consequence, funnel control is feasible for the sine-Gordon equation (\ref{sineGordon}) with the output given in (\ref{OutputSineGordon}).

We shall now illustrate the feasability of funnel control on (\ref{sineGordon}) with some numerical simulations. Let us consider the following set of parameters: $\alpha = \pi+\frac{1}{6}, \nu = -1$. The initial conditions for the variables $x$ and $\partial_tx$ have been chosen as
\begin{align*}
x(0,z) = \frac{1}{6}\sin\left(\frac{\pi}{2}z\right), \partial_t x(0,z) = \frac{1}{5}(2z^2-z^4),
\end{align*} 
while the reference signal $y_{\text{ref}}(t) = \frac{1}{5}\cos(e^{-\frac{t}{4}})$. The function $\phi(t)$, whose inverse determines the funnel boundaries, is fixed to $\frac{1}{e^{-2t}+2.5*10^{-3}}$.

\begin{figure}[h!]
\begin{center}
\includegraphics[scale=0.45,trim=1.5cm 2.8cm 0cm 0.1cm]{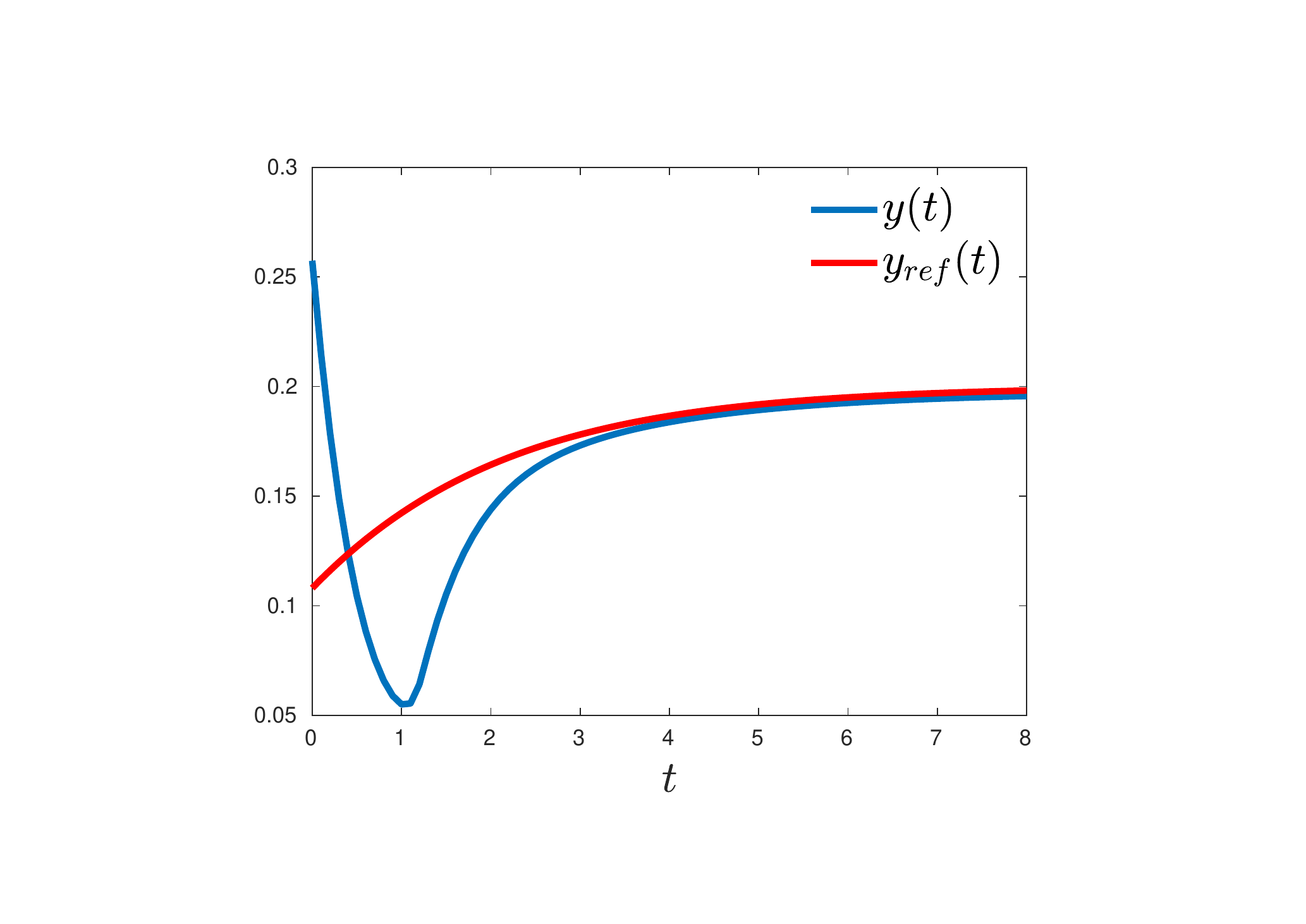}
\end{center}
\caption{Output trajectory (\ref{OutputSineGordon}) with the reference signal $y_{\text{ref}}(t)$.\label{fig:OutputSineGordon}}
\end{figure}

\begin{figure}
\begin{center}
\includegraphics[scale=0.45,trim=1.5cm 2.8cm 0cm 0.1cm]{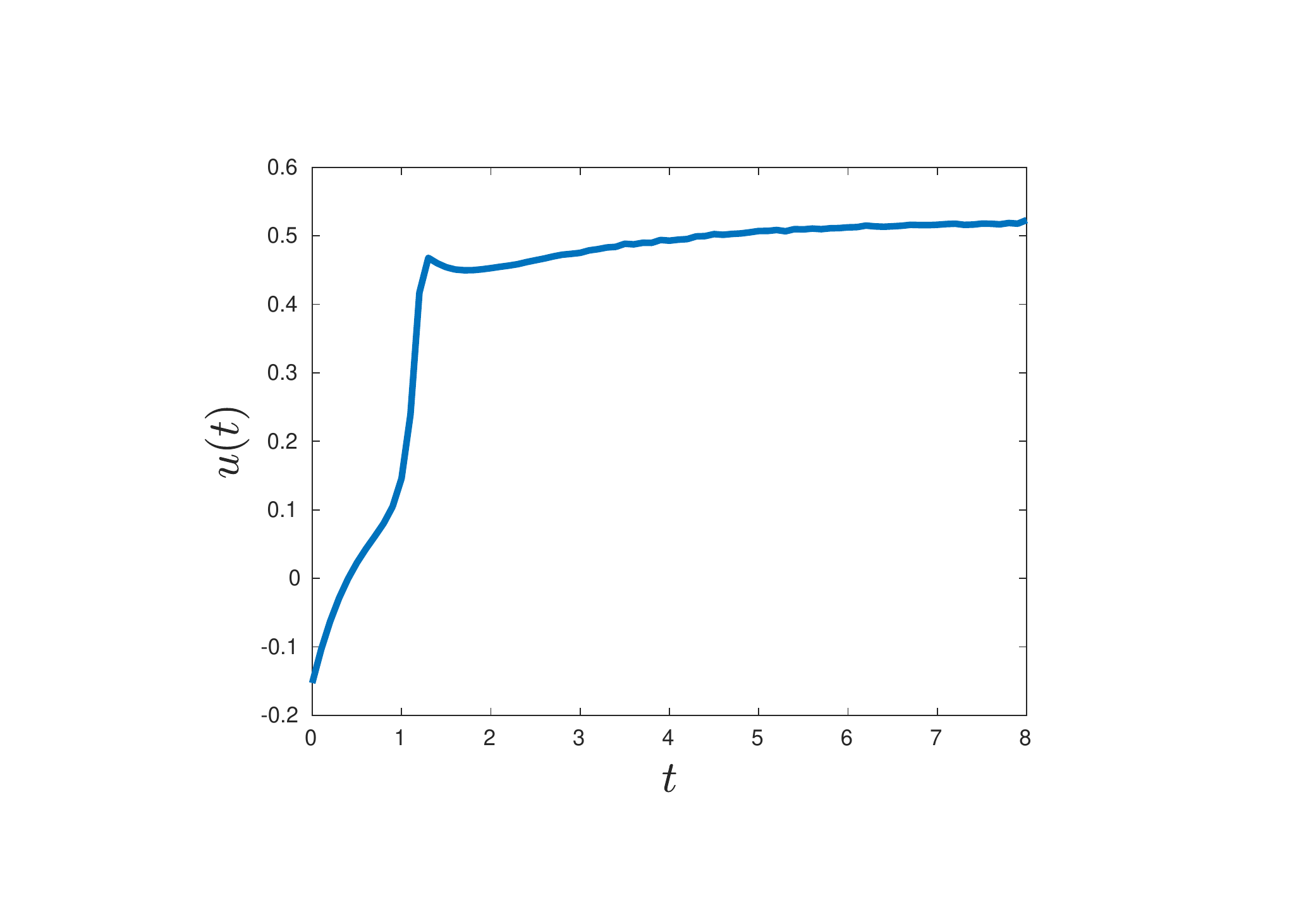}
\end{center}
\caption{Input trajectory.\label{fig:InputSineGordon}}
\end{figure}

\begin{figure}
\begin{center}
\includegraphics[scale=0.45,trim=1.5cm 2.8cm 0cm 0.1cm]{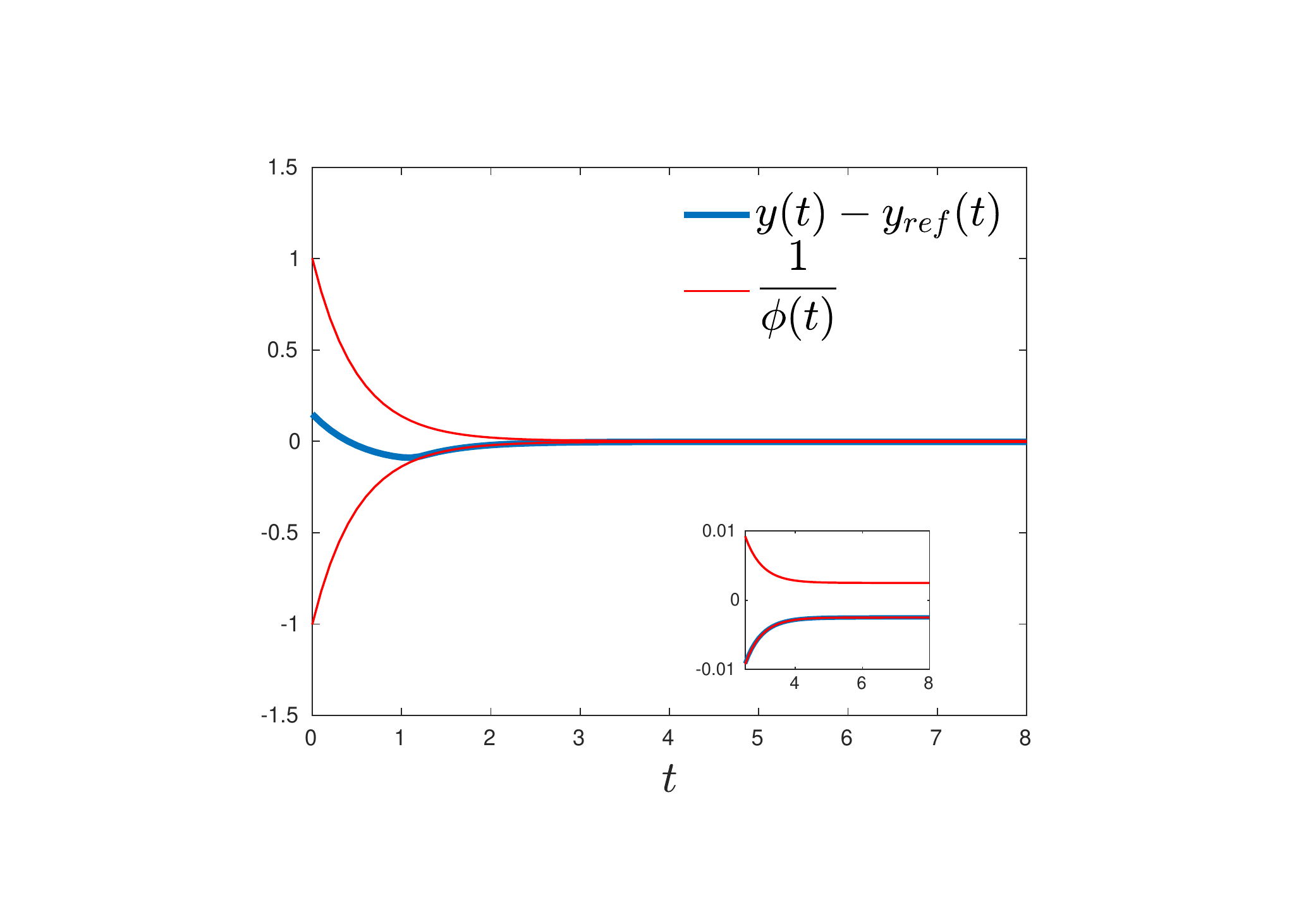}
\end{center}
\caption{Output error trajectory $y(t)-y_{\text{ref}}(t)$ with the funnel boundaries $\frac{1}{\phi(t)}$ and $-\frac{1}{\phi(t)}$.\label{fig:ErrorSineGordon}}
\end{figure}

\begin{figure}
\begin{center}
\includegraphics[scale=0.45,trim=1.5cm 2.8cm 0cm 0.1cm]{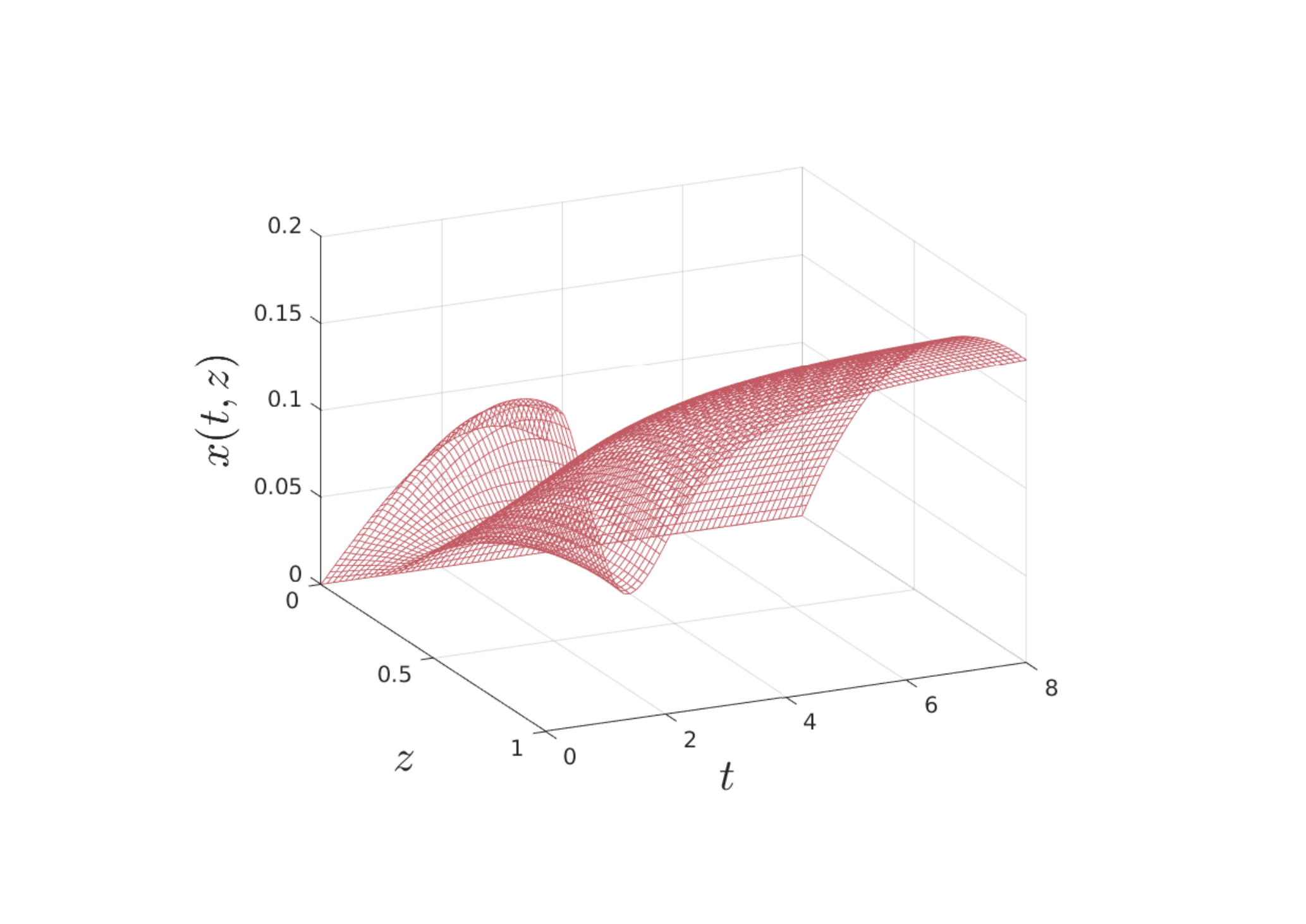}
\end{center}
\caption{Closed-loop state trajectory $x(t,z)$.\label{fig:xSineGordon}}
\end{figure}

\begin{figure}
\begin{center}
\includegraphics[scale=0.45,trim=1.5cm 2.8cm 0cm 0.1cm]{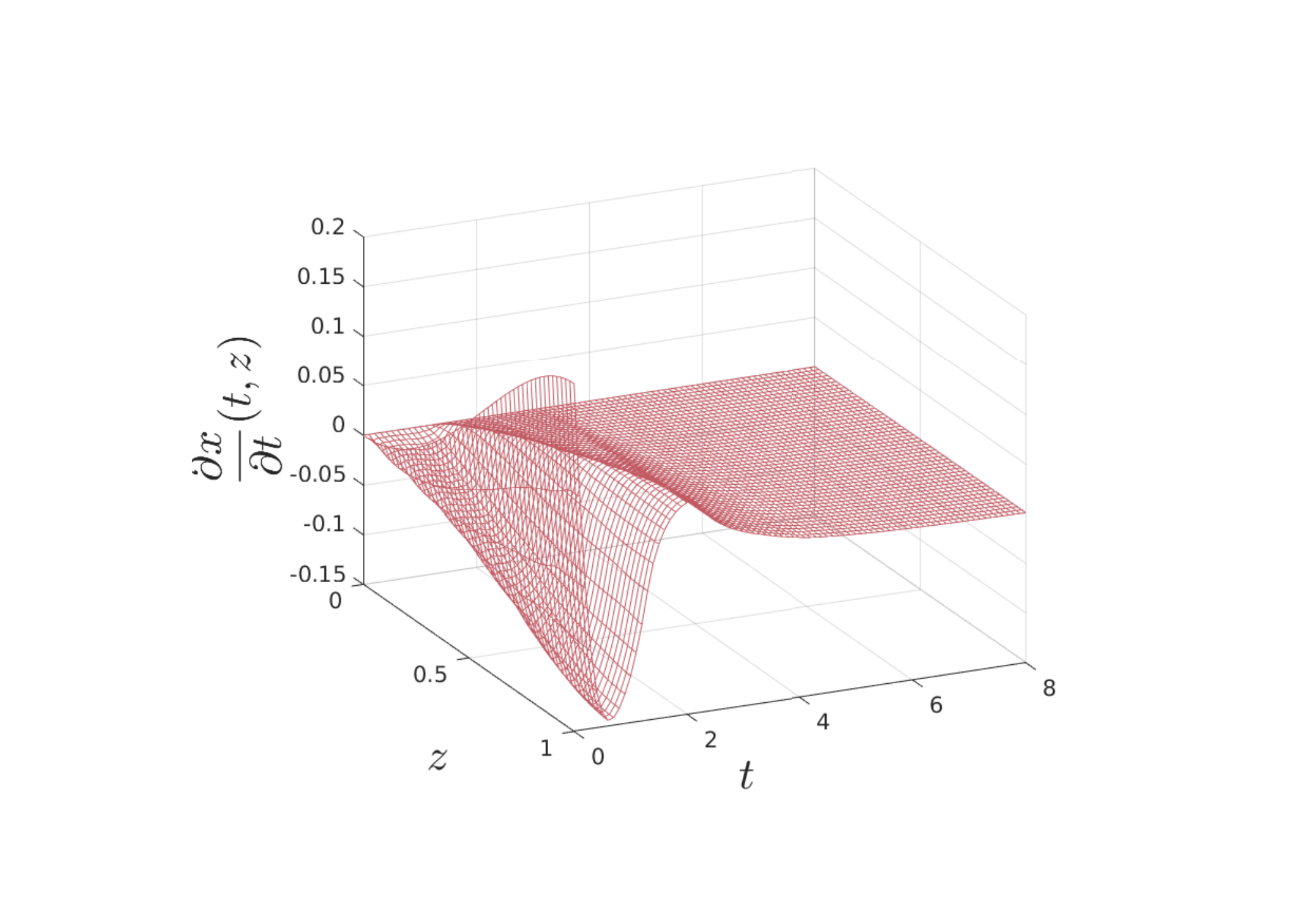}
\end{center}
\caption{Closed-loop state trajectory $\partial_tx(t,z)$.\label{fig:dxSineGordon}}
\end{figure}

The output trajectory (\ref{OutputSineGordon}) with the reference signal $y_{\text{ref}}(t)$ are represented in Figure \ref{fig:OutputSineGordon} whereas the corresponding funnel control is given in Figure \ref{fig:InputSineGordon}. It can be seen that the output tracks the reference quite well. The tracking error is depicted in Figure \ref{fig:ErrorSineGordon} wherein one observes that it remains within the funnel boundaries, as was to be expected.
The state trajectories corresponding to $x(t,z)$ and $\partial_tx(t,z)$ are shown in Figures \ref{fig:xSineGordon} and \ref{fig:dxSineGordon}, respectively.

\section{Conclusion and Perspectives}\label{Conclusions}
Funnel control for a class of nonlinear infinite-dimensional systems was considered. On the basis of the theory presented in \cite{SchwenningerFunnel} and developed in \cite{IlchmannByrnesIsidori} for linear systems with an arbitrary relative degree $r\in\mathbb{N}_0$, we have shown how to extend the feasability of funnel control to the given system class which has to satisfy some well-posedness and BIBO stability assumptions. This constitutes the main result of the paper and is mainly based on a change of variables that enables to extract the output dynamics of the system. This is known as the Byrnes-Isidori form for linear systems, see e.g. \cite{IlchmannByrnesIsidori}. A way of getting the BISBO stability assumption was also established. Two examples whose dynamics are written as a nonlinear PDE and for which the results may be applied were considered. Some numerical simulations reinforce the theoretical results. 

As perspectives for further research, the extension of the results to higher differential orders for the output equation, i.e. in some sense, to systems with higher relative degree could be of great interest. Another perspective aims at weakening the assumption of uniform Lipschitz continuity of the nonlinear operator to only local Lipschitz continuity. As further research the consideration of boundary control and observation could be also an interesting and valuable challenge. This study should include the formalisation of the concept of nonlinear boundary controlled and observed systems. What is meant by well-posedness in this context should also be properly defined and analyzed. For this, we refer to the preliminary works by \cite{Tucsnak}, \cite{Hastir2019SCL} and \cite{Schwenninger2020}. The enlargement of the class of nonlinear infinite-dimensional systems $\Sigma$ for which funnel control is feasible could be interesting to investigate too.

\section*{Acknowledgments}
This research was conducted with the financial support of F.R.S-FNRS. Anthony Hastir is a FNRS Research Fellow under the grant FC 29535. The scientific responsibility rests with its authors.

\section*{References}
\bibliographystyle{elsarticle-harv}       
\bibliography{biblio}   

\end{document}

%% file: FunnelControllerCopy.tex
%\documentclass[tikz, border=2mm]{standalone}

%\usepackage{amsmath}
%\usetikzlibrary{positioning,arrows}
%
%\tikzset{%
%    block/.style={draw, fill=white, rectangle, 
%            minimum height=2em, minimum width=3em},
%    input/.style={inner sep=0pt},       
%    output/.style={inner sep=0pt},      
%    sum/.style = {draw, fill=white, circle, minimum size=2mm, node distance=1.5cm, inner sep=0pt},
%    pinstyle/.style = {pin edge={to-,thin,black}}
%}

\begin{tikzpicture}[auto, node distance=5cm, on grid, >=latex']

\node[input] (input) {};
\node[input, above = of input] (disturbance) {};
\node [block, right = of disturbance] (controller) {$\begin{array}{l}\dot{y}(t) = T(y)(t)+\gamma u(t)+\gamma d(t)\\T(y)(t) = P_0y(t)+S\eta(t)+\langle\tilde{f}(y(t),\eta(t)),c\rangle_H\\y(0)=y_0\end{array}$};
\node [input, above of = controller, node distance=1.5cm] (disturb) {};
\node[input,left = of controller] (input2) {};
\node [input, right = 5cm of controller] (output) {};
\node [block, below = 2.5cm of controller] (funnel) {$\begin{array}{l}
\dot{\eta}(t) = Ry(t) + Q\eta(t) + P^\perp \tilde{f}(y(t),\eta(t))\\
\eta(0) = \eta_0
\end{array}$};
\node [coordinate,right = 4cm of funnel] (eta) {};
\node [coordinate, left = 4cm of controller](test){};
\node [block, below = 2cm of funnel] (Control) {$\frac{-e}{1-\phi^2e^2}$};
\node [sum, right = 5.5cm of Control] (sum1) {};
\node [input, right = 2cm of sum1] (reference) {};
\node [input, right = 2.5cm of Control] (whiteNode) {};
\node [input, above = 0.2cm of whiteNode] (Node) {$e(t)$};

\draw [draw,->] (disturb) node[below right] {$d(t)$} -- (controller);
\draw [draw,-] (controller) -- (output);
\draw [->] (controller) -| node[above right] {$y(t)$} (eta) |- node[left] {} (funnel);
\draw [->] (reference) node[above left,node distance=2cm] {$y_{\text{ref}}(t)$} -- (sum1);
\draw [->] (sum1) -- (Control);
\draw [->] (funnel) -| node {} ([xshift=0mm,yshift=-2mm]test) |-
    node[left] {$\eta(t)$} ([yshift=-2mm]controller.west);
\draw [->] (controller) -| node[] {}
    node [pos=0.99, right] {$-$} (sum1);
\draw [->] (controller) -| node[] {}
    node [pos=0.99, left] {$+$} (sum1);
 \draw [->] (Control) -| node {} ([xshift=-15mm,yshift=2mm]test) |-
    node[left] {$u(t)$} ([xshift=0mm,yshift=2mm]controller.west);

%\node[input] (input) {};
%\node[input, above = of input] (input1) {};
%\node [sum, right = of input] (sum) {};
%\node [block, right = of sum] (controller) {$K(s)$};
%\node [sum, right = of controller] (sum1) {};
%\node [block, right = of sum1] (filterinv) {$H^{-1}(s)$};
%\node [block, right = 2.5cm of filterinv] (system) {$G(s)$};
%\node [output, right = of system] (output) {};
%\node [output, above = of output] (output1) {};
%\node [block, above = of controller] (delay) {$D(s)$};
%\node [sum, below = of sum1] (sum2) {};
%\node [block] (filter) at (sum2-|filterinv) {$H(s)$};
%
%\draw [draw,->] (input) node[above right] {$s_{i-1}$} -- (sum);
%\draw [->] (sum) -- node {$e_{i}$} (controller);
%\draw [->] (controller) -- node {} (sum1);
%\draw [->] (sum1) -- node[name=xi] {$\xi_{i}$} (filterinv);
%\draw [->] (filterinv) -- node[name=u, pos=.3] {$u_{i}$} (system);
%\draw [->] (system) -- (output) node [name=q, above left] {$q_{i}$};
%
%\draw [->] ([xshift=-5mm]q.south) |- (filter);
%\draw [->] (filter) -- node {} (sum2);
%\draw [draw,<-] (sum2) -- ++(90:.6cm) node[above]{$L_i+r_i$};
%
%\draw [->] (sum2) -| node[pos=0.99, right] {$-$} 
%    node [pos=.25, above] {$\tilde{s}_i$} (sum);
%
%\draw [draw,->] (input1) node[above right] (ui-1) {$u_{i-1}$} -- (delay);
%\draw [->] (delay) -| node[] {} 
%    node [near end] {} (sum1);
%
%\draw [->] (u.east|-system) |-  
%    (output1) node[above left] (ui) {$u_i$};
%
%\node[text=red, above left= 5mm and 6mm of ui.west] (veh) {vehicle $i$};
\draw[red, dashed] ([yshift=10mm,xshift=2mm]disturbance.north)-|([xshift=12mm]eta.east)|-([yshift=-10mm]funnel.west)-|([yshift=10mm,xshift=2mm]disturbance.north);%|-(veh.west);
\draw[blue, dashed] ([yshift=8mm,xshift=-5mm]funnel.west)-|([xshift=-1mm]eta.east)|-([yshift=-8mm]funnel.west)-|([yshift=8mm,xshift=-5mm]funnel.west);
\node[text=blue, above left= 10mm and 3mm of funnel] (sigma) {$\Sigma_{y\to\eta}$};
\node[text=red, above left= 13mm and 0mm of output] (sigmaMain) {$\tilde{\Sigma}$};

\end{tikzpicture}